\documentclass[12pt,a4paper]{amsart}

\usepackage[mathscr]{eucal}
\usepackage{amssymb}
\usepackage{amsxtra}
\usepackage{graphics,color,graphicx}
\newcommand{\pA}{\mathcal{A}}
\newcommand{\pB}{\mathcal{B}}

\newcommand{\pL}{\mathcal{L}}

\newcommand{\pS}{\mathcal{S}}

\newcommand{\eC}{\mathscr{C}}
\newcommand{\eF}{\mathscr{F}}
\newcommand{\eI}{\mathscr{I}}
\newcommand{\eK}{\mathscr{K}}
\newcommand{\eL}{\mathscr{L}}
\newcommand{\eM}{\mathscr{M}}
\newcommand{\eN}{\mathscr{N}}
\newcommand{\eH}{\mathscr{H}}
\newcommand{\eU}{\mathscr{U}}
\newcommand{\eX}{\mathscr{X}}

\newcommand{\eY}{\mathscr{Y}}
\newcommand{\eZ}{\mathscr{Z}}
\newcommand{\bC}{\mathbb{C}}
\newcommand{\bR}{\mathbb{R}}

\newcommand{\bN}{\mathbb{N}}

\newcommand{\Ker}{{\rm ker}}
\newcommand{\Image}{{\rm im}}
\newcommand{\Lat}{{\rm Lat}}
\newcommand{\Lath}{{\rm Lat}_{\it h}}
\newcommand{\Latu}{{\rm Lat}_{\it u}}
\newcommand{\Alg}{{\rm Alg}}
\newcommand{\LCom}{{\eC}}
\newcommand{\LInt}{{\eI}}

\newcommand{\IH}{{I_{_{\eH}}}}
\newcommand{\IM}{{I_{_{\eM}}}}
\newcommand{\IN}{{I_{_{\eN}}}}
\newcommand{\IX}{{I_{_{\eX}}}}
\newcommand{\IY}{{I_{_{\eY}}}}
\newtheorem{theorem}{Theorem}
\newtheorem{proposition}[theorem]{Proposition}
\newtheorem{lemma}[theorem]{Lemma}
\newtheorem{corollary}[theorem]{Corollary}
\theoremstyle{definition}
\newtheorem{definition}[theorem]{Definition}
\newtheorem{example}[theorem]{Example}

\begin{document}

%************************************************************************
%************************************************************************
\title[Local commutants and ultrainvariant subspaces]{Local commutants and ultrainvariant subspaces}
\author[J. Bra\v{c}i\v{c}]{Janko Bra\v{c}i\v{c}}
\address{Faculty of Natural Sciences and Engineering, University of Ljubljana, A\v{s}ker\v{c}eva c. 12, SI-1000 Ljubljana, Slovenia}
\email{janko.bracic@ntf.uni-lj.si}
%************************************************************************
\begin{abstract}
For an operator $A$ on a complex Banach space $\eX$ and a closed subspace $\eM\subseteq \eX$, the local commutant of $A$ at $\eM$ is the set $\LCom(A;\eM)$ of all operators $ T$ on $\eX$ such that $TAx=ATx$
for every $x\in \eM$. It is clear that $ \LCom(A;\eM)$ is a closed linear space of operators, however it is not
an algebra, in general. For a given $A$, we show that $\LCom(A;\eM)$ is an algebra if and only if the largest
subspace $\eM_A$ such that $\LCom(A;\eM)=\LCom(A;\eM_A)$ is invariant for every operator in $\LCom(A;\eM)$.
We say that these are ultrainvariant subspaces of $A$. For several types of operators we prove that there exist
non-trivial ultrainvariant subspaces. For a normal operator on a Hilbert space, every hyperinvariant subspace is
ultrainvariant. On the other hand, the lattice of all ultrainvariant subspaces of a non-zero nilpotent operator can be
strictly smaller than the lattice of all hyperinvariant subspaces.
\end{abstract}
%************************************************************************
\maketitle
%************************************************************************

%************************************************************************
%************************************************************************
\section{Introduction} \label{sec01}
%************************************************************************
%************************************************************************

Let $\eX$ be a complex Banach space and let $B(\eX)$ be the Banach algebra of all bounded linear operators on $\eX$.
The dual space of $\eX$ will be denoted by $\eX^*$.
For $A\in B(\eX)$ and a closed subspace $\eM\subseteq \eX$, the local commutant of $A$ at $\eM$ is 
$\LCom(A;\eM)=\{ T\in B(\eX);\; TAx=ATx\;\text{for all}\; x\in \eM\}$. It is clear that $(A)'$, the commutant of $A$, is a subset of $ \LCom(A;\eM)$. Of course, $\LCom(A;\eX)=(A)'$ and $\LCom(A;\{ 0\})=B(\eX)$. 

Local commutants were introduced by Larson \cite{Lar,DL} more than twenty years ago, however there are only a
few results related to them (see \cite{BMZ}, for instance) although that there are some interesting open problems
related to them. One among these problems is the question when is $\LCom(A;\eM)$ an algebra. Of course,
for an arbitrary $A$ and $\eM$ which is equal either to $\{ 0\}$ or $\eX$ the local commutant is an algebra. However, for $A$
which is not a scalar multiple of $\IX$, the identity operator on $\eX$, there exists a closed subspace $\eM\subseteq \eX$ such that
$\LCom(A;\eM)$ is not a subalgebra of $B(\eX)$ but just a linear subspace. 

Recall that a closed subspace $\eM\subseteq \eX$ is said to be invariant for $A\in B(\eX)$ if $A\eM\subseteq \eM$. It is well-known that the family of all closed subspaces which are invariant for $A$ form a lattice; the lattice operations are
the intersection and the closed linear span. As usual, we will denote this lattice by $\Lat(A)$. A sublattice of
$\Lat(A)$ is $\Lath(A)$, the lattice of all hyperinvariant subspaces of $A$: a subspace $\eM \in \Lat(A)$ belongs to
$\Lath(A)$ if and only if $T\eM \subseteq \eM$ for every $T\in (A)'$. It is obvious that the trivial lattice $\{ \{ 0\}, \eX\}$ is a sublattice of $\Lath(A)$ and therefore of $\Lat(A)$. In 1970s, Enflo showed that there exists a Banach
space $\eX$ and a bounded operator $A$ on it such that $\Lat(A)$ is trivial. The result was published much later in \cite{Enf} and meanwhile Read published his example of an operator without a non-trivial invariant subspace, see
\cite{Rea}. Nowadays several examples of Banach spaces (including $l_1$) with operators whose lattice of invariant subspaces is trivial  are known. On the other hand, Lomonosov \cite{Lom} proved that an operator has a non-trivial invariant subspace if it commutes with an operator which is not a scalar multiple of $\IX$ and commutes with a non-zero compact operator. It follows that for some
infinite dimensional Banach spaces $\eX$ (see example given by Argyros and Haydon \cite{AH}) every operator in $B(\eX)$ has a non-trivial invariant subspaces.

However for many Banach spaces (for instance, for the infinite dimensional separable Hilbert space) it is not known if $\Lat(A)$ is non-trivial for every operator $A$. Similarly, for more or less the same class of Banach spaces,
it is not known if $\Lath(A)$ can be trivial for an operator $A$ which is not a scalar multiple of $\IX$. Work on these
problems is still very active, see \cite{Tca} and references therein. The reader is referred to \cite{CP}
and \cite{RR}
for more details about the problem of invariant subspaces.

It turns out that our problem about the algebraic structure of $\LCom(A;\eM)$ is connected with hyperinvariant subspaces of $A$. Our main observation is the following. 
Let  $A\in B(\eX)$ and let $\eM \subseteq \eX$ be a closed subspace. Denote by $\eM_A\subseteq \eX$ the set of all vectors $x\in \eX$ such that $TAx=ATx$ 
for every $T\in \LCom(A;\eM)$. Then $\eM_A$ is a closed subspace of $\eX$ such that $\eM\subseteq \eM_A$ and $\LCom(A;\eM)=\LCom(A;\eM_A)$. The local commutant 
$\LCom(A;\eM)$ is an algebra if and only if $\eM_A$ is invariant for every operator in $\LCom(A;\eM)$; we say that $\eM_A$ is an ultrainvariant subspace of $A$. 
Since $(A)'\subseteq \LCom(A;\eM)$ an ultrainvariant subspace is a hyperinvariant subspace, as well.
The opposite is not true, in general.

Ultrainvariant subspaces are the main theme of the second part of this paper. Denote by $\Latu(A)$ the family
of all ultrainvariant subspaces of $A$. It is a sublattice of $\Lath(A)$. For several types of operators $A$
we show that $\Latu(A)$ is not trivial. There are classes of operators for which
we are able to characterize the lattice of ultrainvariant subspaces completely. For instance, for algebraic operators, in particular, for operators on a finite dimensional Banach space, and for normal
operators on a Hilbert space. If $A$ is a unicellular operator on a Hilbert space, then $\Latu(A)=\Lath(A)=\Latu(A)$.
This holds, for instance, for the Volterra operator and Donoghue operators. 

Several questions about ultrainvariant subspaces remain open. For instance, we do not know if every compact 
operator has a non-trivial ultrainvariant subspace. Actually we believe that it is hard to find an operator which
has a non-trivial hyperinvariant subspace but does not have a non-trivial ultrainvariant subspace.

In the end of this section, let us say a few words about the terminology. Sz.-Nagy and Foia\c{s} initiated the study
of hyperinvariant subspaces in their book {\em Analyse harmonique des op\'{e}rateurs de l'espace Hilbert}
however they called them ultrainvariant subspaces. In this paper we use their terminology (see Definition \ref{def01})
for a special type of hyperinvariant subspaces. Another new term which we use is  {\em girder}. The girder
of $\LCom(A;\eM)$ is the largest subspace $\eM_A\subseteq \eX$ such that $\LCom(A;\eM)=\LCom(A;\eM_A)$.
In our opinion the term girder is suitable as $\eM_A$ ``carries'' the whole space $\LCom(A;\eM)$ in the sense that
$\eM_A$ is the largest subspace whose vectors give the information whether $T$ belongs to $\LCom(A;\eM)$ or not.

%************************************************************************
%************************************************************************
\section{Local commutants} \label{sec02}
\setcounter{theorem}{0}
%************************************************************************
%************************************************************************

Although that we are mainly interested in local commutants we begin our study by considering a slightly more general objects.
Let $\eX$ and $\eY$ be complex Banach spaces and let $B(\eX,\eY)$ be the space of all bounded linear operators from $\eX$ to $\eY$.
For $A\in B(\eX)$ and $B\in B(\eY)$, let $\LInt(A,B)$ be the set of all operators $S\in B(\eX,\eY)$ which intertwine $A$ and $B$, that is,
$\LInt(A,B)=\{ S\in B(\eX,\eY);\; SA=BS\}$. Of course, if $\eX=\eY$ and $A=B$, then $\LInt(A,A)=(A)'$.
Let $\eM \subseteq \eX$ be a closed subspace. The space of operators which intertwine $A$ and $B$ locally at $\eM$ is 
$\LInt(A,B;\eM)=\{ S\in B(\eX,\eY);\; SAx=BSx\; \text{for every}\; x\in \eM\}$. In particular, $\LInt(A,A;\eM)=\LCom(A;\eM)$ is
the local commutant of $A$ at $\eM$. 

It is easily seen that $\LInt(A,B;\eM)$ is closed in the strong operator topology and that
$\LInt(A,B;\eM)=\LInt(A-\lambda \IX,B-\lambda \IY;\eM)$, for every $\lambda \in \bC$. 
For arbitrary invertible operators $U\in B(\eX)$ and $V\in B(\eY)$, we have
\begin{equation} \label{eq06}
\LInt(A,B;\eM)=V^{-1}\LInt(UAU^{-1},VBV^{-1};U\eM)U.
\end{equation}
If $\eM_1, \eM_2$ are closed subspaces of $\eX$ such that $\eM_2 \subseteq \eM_1$, then $\LInt(A,B;\eM_1)\subseteq \LInt(A,B;\eM_2)$. In particular,  
$\LInt(A,B)=\LInt(A,B;\eX)\subseteq \LInt(A,B;\eM)\subseteq \LInt(A,B;\{ 0\})=B(\eX,\eY)$.

Let $\{ \eM_j;\; j\in J\}$ be a non-empty family of closed subspaces of $\eX$.  It is not hard to see that 
$$\bigvee\limits_{j\in J}\{\LInt(A,B;\eM_j);\; j\in J\}\subseteq \LInt(A,B;\bigcap\limits_{j\in J}\eM_j).$$
On the other hand,
\begin{equation} \label{eq05}
\LInt(A,B;\bigvee\limits_{j\in J}\eM_j)=\bigcap\limits_{j\in J}\LInt(A,B;\eM_j).
\end{equation}
Here we denoted by $\bigvee\limits_{j\in J}$ the closed linear span of the involved sets.

\noindent
To verify \eqref{eq05}, note that $\LInt(A,B;\bigvee\limits_{j\in J}\eM_j)\subseteq \bigcap\limits_{j\in J}\LInt(A,B;\eM_j)$
holds since $\eM_i\subseteq \bigvee\limits_{j\in J}\eM_j$ and therefore $\LInt(A,B;\bigvee\limits_{j\in J}\eM_j)\subseteq \LInt(A,B;\eM_i)$,
for every $i\in J$. For the opposite inclusion, observe that $BSx_j=SAx_j$ for every $x_j\in \eM_j$ and every $j\in J$
if $S\in \bigcap\limits_{j\in J}\LInt(A,B;\eM_j)$. Hence $S\in \LInt(A,B;\bigvee\limits_{j\in J}\eM_j)$. (Here we used
a simple fact that $BSx=SAx$ holds for every $x\in \bigvee \eU$ if it holds for every $x\in \eU\subseteq \eX$.)

Let $A\in B(\eX)$, $B\in B(\eY)$ and let $\eM \subseteq \eX$ be a complemented closed subspace, that is, there exists a closed subspace
$\eN \subseteq \eX$ such that $\eX=\eM\oplus \eN$. Let $\left[\begin{smallmatrix}
A_{11} & A_{12}\\ A_{21} & A_{22}
\end{smallmatrix}\right]$ be the operator matrix of $A$ with respect to this decomposition.
For arbitrary operators $S\in B(\eX,\eY)$ and $T\in B(\eX)$, let $\left[ S_1\, S_2\right]$, respectively
$\left[\begin{smallmatrix}
T_{11} & T_{12}\\ T_{21} & T_{22}
\end{smallmatrix}\right]$, be their operator matrices with respect to the decomposition $\eX=\eM \oplus \eN$. These assumptions and notation will be used 
throughout this section.

A straightforward computation shows that $S\in \LInt(A,B)$ if and only if
\begin{equation} \label{eq02}
BS_1=S_1 A_{11}+S_2A_{21}
\end{equation}
and
\begin{equation} \label{eq08}
BS_2=S_1 A_{12}+S_2A_{22}.
\end{equation}
Similarly, $T\in (A)'$ if and only if
\begin{equation} \label{eq09}
T_{11}A_{11}+T_{12}A_{21}=A_{11}T_{11}+A_{12}T_{21},\qquad T_{21}A_{11}+T_{22}A_{21}=A_{21}T_{11}+A_{22}T_{21}
\end{equation}
and
\begin{equation} \label{eq10}
T_{11}A_{12}+T_{12}A_{22}=A_{11}T_{12}+A_{12}T_{22},\qquad T_{21}A_{12}+T_{22}A_{22}=A_{21}T_{12}+A_{22}T_{22}.
\end{equation}

\begin{lemma} \label{lem02}
Operator $S\in B(\eX,\eY)$ is in $\LInt(A,B;\eM)$ if and only if \eqref{eq02} holds and $T\in B(\eX)$ is in $\LCom(A;\eM)$ if and only if
\eqref{eq09} holds.
\end{lemma}

\begin{proof} 
Let $x\in \eX$ be arbitrary and let $x=x_1\oplus x_2$, where $x_1\in \eM$ and $x_2\in \eN$. Then $SAx=(S_1A_{11}+S_2A_{21})x_1+(S_1A_{12}+S_2A_{22})x_2$ and
$BSx=BS_1 x_1+BS_2 x_2$. Hence, if $S\in \LInt(A,B;\eM)$, that is, $BS(x_1\oplus 0)=SA(x_1\oplus 0)$ for every $x_1\in \eM$,
then $BS_1 x_1=(S_1A_{11}+S_2A_{21})x_1$ for every $x_1\in \eM$ which gives \eqref{eq02}. It is clear that
the opposite implication holds, as well. The second part of the assertion follows by the first one if we put $B=\left[\begin{smallmatrix}
A_{11} & A_{12}\\ A_{21} & A_{22}
\end{smallmatrix}\right]$, $S_1=\left[\begin{smallmatrix} T_{11}\\ T_{21}\end{smallmatrix}\right]$ and
$S_2=\left[\begin{smallmatrix} T_{12}\\ T_{22}\end{smallmatrix}\right]$ in \eqref{eq02}.
\end{proof}

It is clear that $\LInt(A,B;\{ 0\})=B(\eX,\eY)$ for arbitrary operators $A$ and $B$. Let $\eM\ne \{0\}$. The following corollary characterizes those pairs of operators $A\in B(\eX)$
and $B\in B(\eY)$ for which $\LInt(A,B;\eM)=B(\eX,\eY)$.

\begin{proposition} \label{prop21}
Let $\eM \ne \{ 0\}$. Then $\LInt(A,B;\eM)=B(\eX,\eY)$ if and only if $A_{11}=\lambda \IM$, $B=\lambda \IY$, for some $\lambda \in \bC$, and $A_{21}=0$.
In particular, $\LCom(A;\eM)=B(\eX)$ if and only if $A$ is a scalar multiple of $\IX$.
\end{proposition}

\begin{proof}
Assume that $\LInt(A,B;\eM)=B(\eX,\eY)$. Let $S=\left[ S_1\, 0\right]$, where $S_1\in B(\eM,\eY)$ is arbitrary. By Lemma \ref{lem02}, $BS_1=S_1A_{11}$. 
Let $y\in \eY$ and $0\ne \xi\in \eM^*$ be arbitrary. Then $y\otimes \xi\in B(\eM,\eY)$ and therefore $B(y\otimes \xi)=(y\otimes \xi)A_{11}$. If $x\in \eM$ is
such that $\langle x,\xi\rangle=1$, then $By=B(y\otimes \xi)x=(y\otimes \xi)A_{11}x=\langle A_{11}x,\xi\rangle y$. Hence, $B\in B(\eY)$ is an operator such that every $0\ne y\in \eY$ is its eigenvector. It follows that there exists $\lambda \in \bC$ such that $B=\lambda \IY$.
Now we have $\lambda S_1=S_1 A_{11}$ for every $S_1\in B(\eM,\eY)$ which gives $A_{11}=\lambda \IM$. Since every $S=\left[ S_1\, S_2\right]\in B(\eX,\eY)$
satisfies \eqref{eq02} we have $S_2A_{21}=0$ for every $S_2\in B(\eN,\eY)$. Hence, $A_{21}=0$. The opposite implication is simple: if $A_{11}=\lambda \IM$,
$B=\lambda \IY$ and $A_{21}=0$, then \eqref{eq02} holds for every $S=\left[S_1\, S_2\right]$.

It follows from the first part of this corollary that in the case when $\eX=\eY$ and $A=B$ we have $\LCom(A;\eM)=B(\eX)$ if and only if $A=\lambda \IX$ for some
$\lambda \in \bC$.
\end{proof}

Equality $\LInt(A,B;\eM)=\LInt(A,B)$ holds if and only if every pair of operators $S_1\in B(\eM,\eY)$ and $S_2\in B(\eN,\eY)$ which satisfies \eqref{eq02} satisfies \eqref{eq08}, as well. In the following proposition we determine
those pairs $A\in B(\eX)$ and $B\in B(\eY)$ for which $\LInt(A,B;\eM)=\LInt(A,B)$ holds for a given subspace 
 $\eM\in \Lat(A)$, $\eM\ne \eX$. We will use the following notation. For $S\in B(\eX, \eY)$, the kernel
of $S$ is denoted by $\Ker(S)$ and the range (image) of $S$ by $\Image(S)$. If $\eZ$ is a non-empty set of vectors
in $\eX$, then by $\eZ^\perp$ we denote its annihilator in $\eX^*$, that is, $\eZ^\perp=\{ \xi\in \eX^*;\; \langle x,\xi\rangle =0\; \text{for all}\; x\in \eZ\}$.

\begin{proposition} \label{prop22}
Let $\eM \ne \eX$. Then $\LInt(A,B;\eM)=\LInt(A,B)$ and $\eM \in \Lat(A)$ if and only if there exists $\lambda \in \bC$ such that
$A=\left[\begin{smallmatrix}
A_{11} & A_{12}\\ 0 & \lambda \IN
\end{smallmatrix}\right]$, $B=\lambda \IY$ and $\overline{\Image(A_{12})}\subseteq \overline{\Image(\lambda \IM-A_{11})}$.
In this case, $\LInt(A,B;\eM)=\{ \left[S_1\,S_2\right]\in B(\eX,\eY);\; \overline{\Image(\lambda \IM-A_{11})}\subseteq \Ker(S_1),\; S_2\in B(\eN,\eY)\}$.
\end{proposition}

\begin{proof}
Assume that $\LInt(A,B;\eM)=\LInt(A,B)$, where $\eX\ne \eM\in \Lat(A)$. Then the operator matrix of $A$ with respect to the decomposition $\eX=\eM\oplus \eN$
is $A=\left[\begin{smallmatrix}
A_{11} & A_{12}\\ 0 & A_{22}
\end{smallmatrix}\right]$. Hence, by Lemma \ref{lem02}, an operator $S=\left[S_1\, S_2\right]\in B(\eX,\eY)$ intertwines $A$ and $B$ at $\eM$ if and only if it
satisfies \eqref{eq02} which in our case means
\begin{equation} \label{eq11}
BS_1=S_1A_{11}.
\end{equation}
It follows that $\left[0\, S_2\right]\in \LInt(A,B;\eM)$, and therefore $\left[0\, S_2\right]\in \LInt(A,B)$, for every $S_2\in B(\eN,\eY)$. Thus, 
$\left[0\,S_2\right]$ satisfies
\eqref{eq08}, that is, $BS_2=S_2A_{22}$ holds for every $S_2\in B(\eN,\eY)$. By Proposition \ref{prop21}, there exists $\lambda \in \bC$ such that
$A_{22}=\lambda \IN$ and $B=\lambda \IY$. Now \eqref{eq11} reads as
\begin{equation} \label{eq12}
S_1(\lambda \IM-A_{11})=0.
\end{equation}
Hence, if $\left[S_1\,S_2\right]$ is an operator in $\LInt(A,\lambda \IY;\eM)$, then $\overline{\Image(\lambda \IM-A_{11})}\subseteq \Ker(S_1)$
and $S_2$ is arbitrary. Note that $\left[S_1\,S_2\right]$ is in $\LInt(A,\lambda\IY)$ if and only if it satisfies \eqref{eq08} which reads as
\begin{equation}\label{eq13}
S_1A_{11}=0.
\end{equation}
in our case. Since $\LInt(A,\lambda \IY;\eM)=\LInt(A,\lambda \IY)$ every $S_1$ which satisfies \eqref{eq12} has to satisfy \eqref{eq13}, as well.
Suppose towards a contradiction that there exists $x=A_{12}u\in \Image(A_{12})$ such that $x\not\in \overline{\Image(\lambda \IM-A_{11})}$.
Let $\xi \in \eM^*$ be such that $\xi\in \overline{\Image(\lambda \IM-A_{11})}^\perp$ and $\langle x,\xi\rangle=1$.
For any $0\ne y\in \eY$,
the rank one operator $y\otimes \xi\in B(\eM,\eY)$ satisfies $(y\otimes \xi)(\lambda \IM-A_{11})=0$, which means $y\otimes \xi \in \LInt(A,\lambda\IY;\eM)$,
however $(y\otimes \xi)A_{12}\ne 0$ since $(y\otimes \xi)A_{12}u=\langle x,\xi\rangle y\ne 0$. This proves that
$\overline{\Image(A_{12})}\subseteq \overline{\Image(\lambda \IM-A_{11})}$.

The opposite implication is easily checked and one gets that $\LInt(A,B;\eM)=\{ \left[S_1\,S_2\right]\in B(\eX,\eY);\; \overline{\Image(\lambda \IM-A_{11})}\subseteq \Ker(S_1),\; S_2\in B(\eN,\eY)\}$.
\end{proof}

\begin{corollary} \label{cor13}
If $\eM \in \Lat(A)$, $\eM \ne \eX$, and $A$ is not a scalar multiple of $\IX$, then $(A)'\subsetneq \LCom(A;\eM)$.
\end{corollary}

\begin{theorem} \label{theo07}
Let $A\in B(\eX)$ and assume that $A$ is not a scalar multiple of $\IX$. Let $\eM\in \Lat(A)$, $\{ 0\}\ne \eM \ne \eX$. The local commutant $\LCom(A;\eM)$
is an algebra if and only if $\LCom(A;\eM)=\{T\in B(\eX);\; T_{11}\in (A_{11})'\; \text{and}\; T_{21}=0\}$.
\end{theorem}

\begin{proof}
We begin the proof with the following simple observation.\\ 
{\em Claim.} Let $\eX_1, \eX_2,\eX_3$ and $\eX_4$ be Banach spaces. If $T\in B(\eX_1,\eX_2)$ and $R\in B(\eX_3,\eX_4)$ are such that $R(u\otimes \xi)T=0$ 
for all $u\in \eX_3$ and $\xi\in \eX_{2}^{*}$, then either $R=0$ or $T=0$.

\noindent
Indeed! Assume that $T\ne 0$. Then there exists $x\in \eX_1$ such that $Tx\ne 0$. Let $\xi\in \eX_{2}^{*}$ be such that $\langle Tx,\xi\rangle=1$.
For an arbitrary $u\in \eX_3$ we have, by the assumption, $R(u\otimes \xi)T=0$ and therefore $R(u\otimes \xi)Tx=Ru=0$.
This shows that $R=0$.

Assume that $\LCom(A;\eM)$ is an algebra. Since $\eM\in \Lat(A)$ the operator matrix of $A$ is of the form 
$\left[\begin{smallmatrix} A_{11} & A_{12}\\ 0 & A_{22}\end{smallmatrix}\right]$. Hence, by Lemma \ref{lem02}, an operator
$T=\left[\begin{smallmatrix} T_{11} & T_{12}\\ T_{21} & T_{22}\end{smallmatrix}\right]$ is in $\LCom(A;\eM)$
if and only if
\begin{equation} \label{eq14}
T_{11}A_{11}-A_{11}T_{11}=A_{12}T_{21},\qquad T_{21}A_{11}-A_{22}T_{21}=0 
\end{equation}
and $T_{12}\in B(\eN,\eM)$, $T_{22}\in B(\eN)$ are arbitrary. In particular, the projection $P$ onto $\eM$ along $\eN$
is in $\LCom(A;\eM)$. Since $\LCom(A;\eM)$ is assumed to be an algebra we see that operators $PTP$ and $(\IX-P)TP$ are in $\LCom(A;\eM)$ for every $T\in \LCom(A;\eM)$. Note that the operator matrices of $PTP$ and $(\IX-P)TP$ are
$\left[\begin{smallmatrix} T_{11} & 0\\ 0 & 0\end{smallmatrix}\right]$ and $\left[\begin{smallmatrix} 0 & 0\\ T_{21} & 0\end{smallmatrix}\right]$,
respectively. Hence, if $T\in \LCom(A;\eM)$, then the entries of its operator matrix have to satisfy
\begin{equation} \label{eq15}
T_{11}A_{11}-A_{11}T_{11}=0,\qquad A_{12}T_{21}=0,\qquad T_{21}A_{11}-A_{22}T_{21}=0.
\end{equation}

Suppose that $A_{12}\ne 0$. Let $u\in \eN$ and $\xi\in \eN^*$ be arbitrary. Then $U$ whose operator matrix is $\left[\begin{smallmatrix} 0 & 0\\ 0 & u\otimes \xi \end{smallmatrix}\right]$
satisfies \eqref{eq14} and therefore $U\in \LCom(A;\eM)$ which gives $U(\IX-P)TP\in \LCom(A;\eM)$. The operator matrix of $U(\IX-P)TP$ is
$\left[\begin{smallmatrix} 0 & 0\\  (u\otimes \xi)T_{21}& 0 \end{smallmatrix}\right]$ and the entries of this matrix satisfy \eqref{eq15}; in particular,
$A_{12}(u\otimes \xi)T_{21}=0$. Since $u\otimes \xi \in B(\eN)$ is an arbitrary operator of rank at most $1$ and $A_{12}\ne 0$ we conclude, by Claim,
that $T_{21}=0$.

Now we may assume that $A_{12}=0$. If $A_{11}=\lambda \IM$ and $A_{22}=\mu \IN$, then $\lambda\ne \mu$ as $A$ is not a scalar multiple of $\IX$. It follows,
by the third equality in \eqref{eq15}, that $T_{21}=0$. Thus, it remains to consider the case when either $A_{11}$ is not a scalar multiple of $\IM$ or
$A_{22}$ is not a scalar multiple of $\IN$. Suppose that this holds. Let $Q\in B(\eN,\eM)$ be arbitrary and let $V\in B(\eX)$ be the operator whose operator
matrix is $\left[\begin{smallmatrix} 0 & Q\\ 0 & 0\end{smallmatrix}\right]$. Then $V\in \LCom(A;\eM)$ and therefore $V(\IX-P(TP\in \LCom(A;\eM)$
for every $T\in \LCom(A;\eM)$. The operator matrix of $V(\IX-P(TP\in \LCom(A;\eM)$ is $\left[\begin{smallmatrix} QT_{21} & 0\\ 0 & 0\end{smallmatrix}\right]$. 
Hence, $QT_{21}\in (A_{11})'$. Since $T_{21}$ satisfies the third equality in \eqref{eq15} we have $A_{11}QT_{21}=QT_{21}A_{11}=QA_{22}T_{21}$ which
gives $(A_{11}Q-QA_{22})T_{21}=0$. We have supposed that either $A_{11}$ is not a scalar multiple of $\IM$ or $A_{22}$ is not a scalar multiple of $\IN$.
Hence, by Proposition \ref{prop21}, there exists $Q_0\in B(\eN,\eM)$ such that $R=A_{11}Q_0-Q_0A_{22}\ne 0$. Thus, $RT_{21}=0$ for every $T\in \LCom(A;\eM)$.
Let $T$ be fixed now. As in the previous paragraph, let $u\in \eN$ and $\xi\in \eN^*$ be arbitrary and let $U\in \LCom(A;\eM)$ be the operator
whose operator matrix is $\left[\begin{smallmatrix} 0 & 0\\ 0 & u\otimes \xi \end{smallmatrix}\right]$. Then $U(\IX-P)TP\in \LCom(A;\eM)$ and therefore
$R(u\otimes \xi)T_{21}=0$. Since $u\otimes \xi \in B(\eN)$ is an arbitrary operator of rank at most $1$ and $R\ne 0$ we conclude, by Claim,
that $T_{21}=0$. This proves that $\LCom(A;\eM)=\{T\in B(\eX);\; T_{11}\in (A_{11})'\; \text{and}\; T_{21}=0\}$ if $\LCom(A;\eM)$ is assumed to be an algebra.
The opposite implication is obvious.
\end{proof}

In the proof of Theorem \ref{theo07} we have seen that the projection $P$ which maps onto $\eM$ along $\eN$ is in $\LCom(A;\eM)$ if $\eM\in \Lat(A)$.
It is not hard to see that the opposite implication holds as well.

\begin{corollary} \label{cor14}
Projection $P$ is in $\LCom(A;\eM)$ if and only if $\eM$ is invariant for $A$.
\end{corollary}

Let $\sigma(T)$ denote the spectrum of $T\in B(\eX)$. It follows from a well-known theorem proved by Rosenblum that $\LInt(A,B)=\{ 0\}$ if 
$\sigma(A)\cap \sigma(B)=\emptyset$, see \cite[Corollary 3.3]{Ros} and \cite[Theorem I.4.1]{GGK1}. Note that the condition $\sigma(A)\cap \sigma(B)=\emptyset$
is not necessary for $\LInt(A,B)=\{ 0\}$. For instance, let $B=0$ and let $A\in B(\eX)$ be such that $\overline{\Image(A)}=\eX$ and $\Ker(A)\ne \{ 0\}$. Then
$0\in \sigma(A)\cap \sigma(B)$, however $\LInt(A,0)=\{ 0\}$ because $SA=0$ holds for $S\in B(\eX,\eY)$ if and only if $S=0$ as $\overline{\Image(A)}=\eX$. 
On the other hand, note that $\LInt(0,A)\ne \{ 0\}$. Indeed, if $S\in B(\eY,\eX)$ is such that $\Image(S)\subseteq \Ker(A)$, then $AS=0$.

\begin{corollary} \label{cor05}
If $\eM\in \Lat(A)$ and $\sigma(A_{11})\cap \sigma(A_{22})=\emptyset$, then $\LCom(A;\eM)$ is an algebra.
\end{corollary}

\begin{proof}
By Lemma \ref{lem02}, $T\in \LCom(A;\eM)$ if and only if \eqref{eq09} holds. Since $A_{21}=0$ the equalities in \eqref{eq09} simplify to \eqref{eq14}.
Because of $\sigma(A_{11})\cap \sigma(A_{22})=\emptyset$, only $T_{21}=0$ satisfies the second equation in \eqref{eq14}. Thus, 
$\LCom(A;\eM)=\{T\in B(\eX);\; T_{11}\in (A_{11})'\; \text{and}\; T_{21}=0\}$ and therefore it is an algebra, by Theorem \ref{theo07}.
\end{proof}

In the following example we will see that $\LCom(Q;\eL)$ can be an algebra for $Q\in B(\eX)$ and a closed subspace $\eL\subseteq \eX$ which is not 
invariant for $Q$. 

\begin{example} \label{ex01}
Let $\eX_1, \eX_2, \eX_3$ be closed subspaces such that $\eX_1, \eX_2$ are non-trivial. Let $\eX=\eX_1\oplus \eX_2\oplus \eX_3$ and suppose that there exists a surjective 
operator in $B(\eX_1,\eX_2)$. Let $Q\in B(\eX)$ be the projection onto $\eX_1$ along $\eX_2\oplus \eX_3$. Choose and fix a surjective operator $U\in B(\eX_1, \eX_2)$.
Let $\eL \subseteq \eX$ be the subspace of all vectors $x\in \eX$ which are of the form $x=x_1\oplus Ux_1\oplus x_3$, where $x_1\in \eX_1$ and $x_3\in \eX_3$ 
are arbitrary. Denote by $P_j:\eX \to \eX_j $ ($j=1,2,3$) the operators given by $P_j(x_1\oplus x_2\oplus x_3)=x_j$. It is clear that $P_1\eL=\eX_1$ and 
$P_3\eL=\eX_3$. Since $U$ is surjective we have $P_2 \eL=\eX_2$, as well. Note that $\eL \not\in \Lat(Q)$. Namely, if $x=x_1\oplus Ux_1\oplus x_3\in \eL$ is such that
$Ux_1\ne 0$, then $Qx=x_1\oplus 0\oplus 0\not\in \eL$.

Assume that $T\in \LCom(Q;\eL)$. Let $\left[T_{ij}\right]_{i,j=1}^{3}$ be its operator matrix with respect to the decomposition $\eX=\eX_1\oplus \eX_2\oplus \eX_3$. 
For an arbitrary $x=x_1\oplus Ux_1\oplus x_3\in \eL$, we have $TQx=T_{11}x_1 \oplus T_{21}x_1\oplus T_{31}x_1$ and $QTx=(T_{11}x_1+T_{12}Ux_1+T_{13}x_3)\oplus 0\oplus 0$.
Since $x_1\in \eX_1$ and $ x_3\in \eX_3$ are arbitrary and $U$ is surjective it follows from $TQx=QTx$ that $T_{21}=0$, $T_{31}=0$, $T_{12}=0$ and $T_{13}=0$. 
Hence $T\in (Q)'$, that is $\LCom(Q;\eL)=(Q)'$; in particular, $\LCom(Q;\eL)$ is an algebra. 

Note that this example shows that $\LCom(A;\eM)$ and  $(A)'$ can be equal for an operator which is not a scalar multiple of $\IX$ and $\eM \ne \eX$ which is not invariant for $A$ (cf. Corollary \ref{cor13}).~\hfill $\Box$
\end{example}

Although that the condition $\eM\in \Lat(A)$ does not need to be fulfilled for $\LCom(A;\eM)$ being an algebra the invariance is not far away as we shall see in the following section.

%************************************************************************
%************************************************************************
\section{Multiplicative structure of the space of local intertwiners} \label{sec03}
\setcounter{theorem}{0}
%************************************************************************
%************************************************************************

Let $A\in B(\eX)$ and let $\eM\subseteq \eX$ be a closed subspace. We know that $\LCom(A;\eM)$ is not an algebra in general. However it has some non-trivial multiplicative structure.
For instance, let $T\in (A)'$ be such that $\eM\in \Lat(T)$. Then $ ASTx=SATx=STAx$ for all $S\in \LCom(A;\eM)$ and $x\in \eM$. Hence, $\LCom(A;\eM)T\subseteq \LCom(A;\eM)$.
The following proposition is a generalization of this simple observation.

\begin{proposition} \label{prop04}
Let $A\in B(\eX)$, $B\in B(\eY)$ and let $\eM\subseteq \eX$ be a closed subspace.
Inclusion $\LInt(A,B;\overline{T\eM})T\subseteq \LInt(A,B;\eM)$ holds for every $T\in \LCom(A;\eM)$. If $T$ is invertible,
then $\LInt(A,B; \overline{T\eM})T=\LInt(A,B;\eM)$.
\end{proposition}

\begin{proof}
Let $T\in \LCom(A;\eM)$. For arbitrary $S\in \LInt(A,B;\overline{T\eM})$ and $x\in \eM$, we have $STAx=SATx=BSTx$ and therefore $ST\in \LInt(A,B;\eM)$.

Assume now that $T\in \LCom(A;\eM)$ is invertible. Let $S\in \LInt(A,B;\eM)$ be arbitrary. Then for every $x\in \eM$ we have
$ST^{-1}ATx=SAx=BSx=BST^{-1}Tx$ which implies that $ST^{-1}\in \LInt(A,B;\overline{T\eM})$, that is, $S\in \LInt(A,B;\overline{T\eM})T$.
\end{proof}

For a non-empty set of operators $\pS\subseteq B(\eX,\eY)$ and a closed subspace $\eM\subseteq \eX$, let $\pS\eM$ denote the closed linear span
in $\eY$ of all subspaces $S\eM$, where $S\in \pS$, that is, $\pS\eM=\bigvee_{S\in \pS}S\eM.$

\begin{theorem} \label{theo01}
Let $A\in B(\eX)$, $B\in B(\eY)$ and let $\eM \subseteq \eX$ be a closed subspace. Then
$ \LCom\bigl(B;\LInt(A,B;\eM)\eM\bigr)$ is a subalgebra of $B(\eY)$ and it is the largest
subalgebra of $B(\eY)$ such that $\LInt(A,B;\eM)$ is a left module over it.
\end{theorem}

\begin{proof}
Let us show that $\LCom\bigl(B;\LInt(A,B;\eM)\eM\bigr)\cdot \LInt(A,B;\eM)\subseteq \LInt(A,B;\eM)$. 
Assume that $T\in \LCom\bigl(B;\LInt(A,B;\eM)\eM\bigr)$ and $S\in \LInt(A,B;\eM)$. 
For every $x\in \eM$, we have $BSx=SAx$ and therefore $BTSx=TBSx=TSAx$ because $Sx\in \LInt(A,B;\eM)\eM$. We conclude that $TS\in \LInt(A,B;\eM)$.

Let $T_1, T_2\in \LCom\bigl(B;\LInt(A,B;\eM)\eM\bigr)$. If $S\in \LInt(A,B;\eM)$, then $T_2S\in \LInt(A,B;\eM)$, by the inclusion from the previous
paragraph, and, by the same reason, $T_1T_2S\in \LInt(A,B;\eM)$. Hence, for every $S\in \LInt(A,B;\eM)$ and every $x\in \eM$, we have $BT_1T_2Sx=T_1T_2SAx=T_1T_2BSx$. 
Because of linearity and countinuity we may conclude that $T_1T_2$ commutes with $B$ on the subspace $\LInt(A,B;\eM)\eM$ of $\eY$, that is, 
$T_1T_2\in \LCom\bigl(B;\LInt(A,B;\eM)\eM\bigr)$.
Thus, $\LCom\bigl(B;\LInt(A,B;\eM)\eM\bigr)$ is a subalgebra of $\pB(\eY)$ and $\LInt(A,B;\eM)$ is a left module over it.

Assume that $T\in B(\eY)$ is such that $T\LInt(A,B;\eM)\subseteq \LInt(A,B;\eM)$. Let $S\in \LInt(A,B;\eM)$ and $x\in \eM$ be arbitrary. Then
$BTSx=TSAx=TBSx$ which implies $T\in \LCom\bigl(B;\LInt(A,B;\eM)\eM\bigr)$. This shows that
$ \LCom\bigl(B;\LInt(A,B;\eM)\eM\bigr)$ is the largest
subalgebra of $B(\eY)$ such that $\LInt(A,B;\eM)$ is a left module over it.
\end{proof}

\begin{corollary} \label{cor01}
Let $A\in B(\eX)$ and let $\eM\subseteq \eX$ be a closed subspace. Then $\LCom(A;\LCom(A;\eM)\eM)$ is the largest algebra contained
in $\LCom(A;\eM)$. If $\eM\in \Lat\bigl(\LCom(A;\eM)\bigr)$, then $\LCom(A;\eM)$ is an algebra.
\end{corollary}

\begin{proof}
By Theorem \ref{theo01}, $\LCom(A;\LCom(A;\eM)\eM)$ is the largest algebra contained in $\LCom(A;\eM)$. If $\eM\in \Lat\bigl(\LCom(A;\eM)\bigr)$,
then $\LCom(A;\eM)\eM=\eM$ and therefore $\LCom(A;\eM)$ is an algebra.
\end{proof}

Now we consider the right module structure of $\LInt(A,B;\eM)$. 
Let $A\in B(\eX)$, $B\in B(\eY)$ and let $\pS\subseteq B(\eX,\eY)$ be an arbitrary set of operators. It is clear that every $S\in \pS$ intertwines $A$ and $B$
at vector $0$, that is, $\pS\subseteq \LInt(A,B;\{ 0\})$. Hence, the family $\eF$ of all closed subspaces $\eM \subseteq \eX$ such that $\pS\subseteq \LInt(A,B;\eM)$
is not empty. Let $\eM_\pS=\bigvee_{\eM\in \eF}\eM$. By \eqref{eq05}, $\pS \subseteq \bigcap_{\eM\in \eF}\LInt(A,B;\eM)=\LInt(A,B;\eM_{\pS})$,
that is, $\eM_\pS\in \eF$. Of course, $\eM_\pS$ is the largest subspace in $\eF$.
In particular, there exists the largest subspace $\eM_{\LInt(A,B;\eM)}$ of $\eX$ on which all operators from $\LInt(A,B;\eM)$ intertwine $A$ and $B$. 
We will simplify the notation to $\eM_{_{A,B}}$ instead of $\eM_{\LInt(A,B;\eM)}$ (and $\eM_{_A}$ instead of $\eM_{_{A,A}}$ if $\eX=\eY$ and $A=B$)
and we will call this subspace a {\em girder}, more precisely, $\eM_{_{A,B}}$ is the girder of $\LInt(A,B;\eM)$ (and $\eM_A$
is the girder of $\LCom(A;\eM)$) induced by $\eM$. Of course, $\eX$ is the girder of $\LInt(A,B)$. On the other hand, $\{ 0\}$ is not necessary a girder --- think about $\IX$.
However, by Proposition \ref{prop21}, $\{ 0\}$ is the girder of $\LCom(A;\{ 0\})=B(\eX)$ if $A$ is not a scalar multiple of $\IX$.

It is clear that $\bigcap_{S\in \LInt(A,B;\eM)}\Ker(BS-SA)$ is a closed subspace of $\eX$ such that $BSx=SAx$ for every $x$ in this subspace.
Hence $\bigcap_{S\in \LInt(A,B;\eM)}\Ker(BS-SA)\subseteq \eM_{_{A,B}}$. On the other hand, if $x\in \eM_{_{A,B}}$, then $BSx=SAx$ for every 
$S\in \LInt(A,B;\eM)$, by the definition of $\eM_{_{A,B}}$. Hence,
\begin{equation} \label{eq04}
\eM_{_{A,B}}=\bigcap_{S\in \LInt(A,B;\eM)}\Ker(SA-BS).
\end{equation}
It follows that
\begin{equation} \label{eq16}
\LInt(A,B;\eM_{_{A,B}})=\LInt(A,B;\eM)\qquad \text{and}\qquad (\eM_{_{A,B}})_{_{A,B}}=\eM_{_{A,B}}.
\end{equation}
Indeed, since $\eM \subseteq \eM_{_{A,B}}$ we have $\LInt(A,B;\eM_{_{A,B}})\subseteq \LInt(A,B;\eM)$. On the other hand, the opposite inclusion follows by 
\eqref{eq04}. Thus, $\LInt(A,B;\eM_{_{A,B}})=\LInt(A,B;\eM)$ and therefore $(\eM_{_{A,B}})_{_{A,B}}=\eM_{_{A,B}}$, by \eqref{eq04}.

\begin{proposition} \label{prop19}
Let $A\in B(\eX)$, $B\in B(\eY)$ and let $\eK, \eM\subseteq \eX$ be closed subspaces. If $\eM \subseteq \eK$, then $\eM_{_{A,B}}\subseteq \eK_{_{A,B}}$.
In particular, if $\eM \subseteq \eK \subseteq \eM_{_{A,B}}$, then $\eM_{_{A,B}}=\eK_{_{A,B}}$ and therefore $\LInt(A,B;\eM_{_{A,B}})=\LInt(A,B;\eK)$.
\end{proposition}

\begin{proof}
Assume that $\eM \subseteq \eK$. Then $\LInt(A,B;\eM)\supseteq \LInt(A,B;\eK)$. Hence, if $x\in \eM_{_{A,B}}$, then $SAx=BSx$ for every
$S\in \LInt(A,B;\eK)$ and therefore $x\in \eK_{_{A,B}}$.

If $\eM \subseteq \eK \subseteq \eM_{_{A,B}}$, then $\LInt(A,B;\eM)\supseteq \LInt(A,B;\eK)\supseteq \LInt(A,B;\eM_{_{A,B}})=\LInt(A,B;\eM)$, that is,
$\LInt(A,B;\eM_{_{A,B}})=\LInt(A,B;\eK)$. It follows from this and the definition of $\eK_{_{A,B}}$ that $\eK_{_{A,B}}=\eM_{_{A,B}}$.
\end{proof}

For a closed subspace $\eM\subseteq \eX$, we denote by $\Alg(\eM)$ the algebra of all operators 
$T\in B(\eX)$ such that $\eM \in \Lat(T)$.

\begin{theorem} \label{theo02}
Let $A\in B(\eX)$, $B\in B(\eY)$ and let $\eM \subseteq \eX$ be a closed subspace. Then
$\LCom(A;\eM_{_{A,B}})\cap \Alg(\eM_{_{A,B}})$ is a subalgebra of $B(\eX)$ and $\LInt(A,B;\eM)$ is a right module over it.

If $\eX=\eY$ and $A=B$, then $\LCom(A;\eM_{_{A}})\cap \Alg(\eM_{_{A}})$ is the largest subalgebra of $B(\eX)$ such that
$\LCom(A;\eM)$ is a right module over it.
\end{theorem}

\begin{proof}
Recall that $\eM \subseteq \eM_{_{A,B}}$, by the definition of $\eM_{_{A,B}}$, and $\LInt(A,B;\eM)=\LInt(A,B;\eM_{_{A,B}})$,
by \eqref{eq16}. Let $S\in \LInt(A,B;\eM)$ and $T\in \LCom(A;\eM_{_{A,B}})\cap \Alg(\eM_{_{A,B}})$ be arbitrary. If $x\in \eM$,
then $STAx=SATx=BSTx$. The first equality holds because of $T\in \LCom(A;\eM_{_{A,B}})$ and the second because of $Tx\in \eM_{_{A,B}}$ and $S\in \LInt(A,B;\eM_{_{A,B}})$. 
We have proved that
$$ \LInt(A,B;\eM)\cdot\bigl(\LCom(A;\eM_{_{A,B}})\cap \Alg(\eM_{_{A,B}})\bigr)\subseteq \LInt(A,B;\eM).$$

If $T_1, T_2\in \LCom(A;\eM_{_{A,B}})\cap \Alg(\eM_{_{A,B}})$, then, of course, $T_1T_2\in \Alg(\eM_{_{A,B}})$. Let $x\in \eM_{_{A,B}}$
be arbitrary. Since $T_2\in \LCom(A;\eM_{_{A,B}})$ we have $T_1T_2Ax=T_1AT_2x$. Because of $T_2x\in \eM_{_{A,B}}$ we have $AT_1T_2x=T_1AT_2x$, as well. 
Hence, $T_1T_2Ax=AT_1T_2x$ for every $x\in \eM_{_{A,B}}$, that is, 
$T_1T_2\in \LCom(A;\eM_{_{A,B}})$. We have proved that $\LCom(A;\eM_{_{A,B}})\cap \Alg(\eM_{_{A,B}})$ is a subalgebra of $B(\eX)$ and
$\LInt(A,B;\eM)$ is a right module over it.

Assume now that $\eX=\eY$ and $A=B$. By the first part of the proof we already know that $\LCom(A;\eM_{_{A}})\cap \Alg(\eM_{_{A}})$
is a subalgebra of $B(\eX)$ and $\LCom(A;\eM_{_A})$ is a right module over it. Let $T\in B(\eX)$ be an arbitrary operator such that
$\LCom(A;\eM_{_A}) T\subseteq \LCom(A;\eM_{_A})$. Since $(A)'\subseteq \LCom(A;\eM_{_A})$ and $(A)'$ contains the identity operator $I$
we see that $T\in \LCom(A;\eM_{_A})$. Let $x\in \eM_{_A}$ and $S\in \LCom(A;\eM_{_A})$ be arbitrary. Then $SATx=STAx$. Because of
$ST\in \LCom(A;\eM_{_A})$ we also have $ASTx=STAx$. Hence $(SA-AS)Tx=0$ for all $S\in \LCom(A;\eM_{_A})$ which means, see \eqref{eq04},
that $Tx\in \eM_{_A}$. This proves that $T\in \Alg(\eM_{_A})$. Thus, $\LCom(A;\eM_{_{A}})\cap \Alg(\eM_{_{A}})$ is the largest subalgebra of $B(\eX)$ such that
$\LCom(A;\eM)$ is a right module over it.
\end{proof}

In the following theorem we list equivalent conditions for $\LCom(A;\eM)$ being an algebra. 

\begin{theorem} \label{theo04}
Let $A\in B(\eX)$, $A\ne \lambda \IX$ for every $\lambda \in \bC$, and let $\eM \subseteq \eX$ be a closed subspace. The following assertions are equivalent:

(i) $\LCom(A;\eM)$ is an algebra;

(ii) $\LCom(A;\eM) \eM \subseteq \eM_A$;

(iii) $\LCom(A;\eM) \eM = \eM_A$;

(iv) $\eM_{_A}\in \Lat\bigl(\LCom(A;\eM)\bigr)$.
\end{theorem}

\begin{proof}
$(i)\Rightarrow(ii)$. Let $S,T\in \LCom(A;\eM)$ be arbitrary. Since $ST\in \LCom(A;\eM)$ we have
$(SA-AS)Tx=STAx-STAx=0$ for every $x\in \eM$. Hence, $T x\in \Ker(SA-AS)$. Since $S\in \LCom(A;\eM)$ is arbitrary we have
$Tx\in \eM_{_A}$ which means that $T\eM \subseteq \eM_{_A}$. Now we may conclude that $\LCom(A;\eM)\eM\subseteq \eM_{_A}$.

$(ii)\Rightarrow(iii)$. If $\LCom(A;\eM)\eM=\eX$, then $\LCom(A;\eM)\eM=\eM_A$, of course. Assume therefore that $\LCom(A;\eM)\eM\ne \eX$.
Then there exist $z\in \eX \setminus \LCom(A;\eM)\eM$ and 
$\xi \in (\LCom(A;\eM)\eM)^\perp$ such that $\langle z,\xi\rangle=1$. Since $A$ is not a scalar multiple of $\IX$ there exists $0\ne y\in \eX$  which is not an eigenvector of $A$. 
Since $Ax\in \LCom(A;\eM)\eM$ whenever $x\in \LCom(A;\eM)\eM$ we have $A(y\otimes \xi)x=0=(y\otimes \xi)Ax$ for every $x\in \LCom(A;\eM)\eM$. Hence,
$y\otimes \xi \in \LCom(A;\eM)$. It follows from $A(y\otimes \xi)z=Ay$ and $(y\otimes \xi)Az=\langle Az,\xi\rangle y$ and the assumption that $y$ is
not an eigenvector of $A$ that $A(y\otimes \xi)z\ne (y\otimes \xi)Az$ which implies $z\not\in \eM_A$.

$(iii)\Rightarrow(i)$.  Let $S, T\in \LCom(A;\eM)$ be arbitrary. Then, for every $x\in \eM$, we have
$(STA-AST)x=SATx-SATx=0$ since it follows by the assumption that $Tx\in \eM_{_A}$ and we already know, see \eqref{eq16}, that 
$\LCom(A;\eM)=\LCom(A;\eM_{_A})$. Hence, $ST\in \LCom(A;\eM)$.

$(i)\Rightarrow(iv)$. Since $\LCom(A;\eM)=\LCom(A;\eM_{_A})$ we have $\LCom(A;\eM_{_A})\eM_{_A}\subseteq \bigl(\eM_{_A}\bigr)_{_A}=\eM_{_A}$,
by the equivalence of assertions (i) and (ii). Hence, $\eM_{_A}\in \Lat\bigl(\LCom(A;\eM)\bigr)$. 

$(iv)\Rightarrow(i)$. If $\eM_{_A}\in \Lat\bigl(\LCom(A;\eM)\bigr)$,
then $\LCom(A;\eM)\eM\subseteq \LCom(A;\eM)\eM_{_A}\subseteq \eM_{_A}$ and therefore $\LCom(A;\eM)$ is an algebra, by equivalence of assertions (i) and (ii).
\end{proof}

%************************************************************************
%************************************************************************
\section{Ultrainvariant subspaces} \label{sec04}
\setcounter{theorem}{0}
%************************************************************************
%************************************************************************

We begin this section with a definition.

\begin{definition} \label{def01}
A closed subspace $\eM \subseteq \eX$ is {\em ultrainvariant} for $A\in B(\eX)$ if it is invariant for every operator in $\LCom(A;\eM)$.
\end{definition}

Hence, every ultrainvariant subspace of $A$ is a hyperinvariant subspace, as well. We will see that the opposite implication does not hold in general. For instance, if $A\ne 0$ is a nilpotent operator, then $\overline{\Image(A^k)}$ is a hyperinvariant subspace
of $A$ for every $k\in \bN$, however it is not necessary an ultrainvariant subspace (see Theorem \ref{theo03}).

\begin{proposition} \label{prop06}
Let $A\in B(\eX)$, $A\ne \lambda \IX$ for every $\lambda \in \bC$. A closed subspace
$\eM\subseteq \eX$ is ultrainvariant for $A$ if and only if $\LCom(A;\eM)$ is an algebra and $\eM$ is its girder.
\end{proposition}

\begin{proof}
Assume that $\eM$ is an ultrainvariant subspace of $A$. Let $T_1, T_2\in \LCom(A;\eM)$ be arbitrary. Then $T_1T_2Ax=T_1AT_2x=AT_1T_2x$ for every $x\in \eM$ because
$T_2x\in \eM$. Thus, $\LCom(A;\eM)$ is an algebra. By Theorem \ref{theo04}, $\LCom(A;\eM)\eM=\eM_A$. However, $\LCom(A;\eM)\eM=\eM$ as $\eM$ is invariant for every operator in $\LCom(A;\eM)$, that is, $\eM=\eM_A$. The opposite implication follows by Theorem \ref{theo04}.
\end{proof}

Subspaces $\{ 0\}$ and $\eX$ are trivial ultrainvariant subspaces of every $A\in B(\eX)$. If $A$
is a scalar multiple of $\IX$, then it is obvious that $A$ has only trivial ultrainvariant subspaces as it has only trivial hyperinvariant subspaces. 

In this section we will prove a few general results about ultrainvariant subspaces. Examples of ultrainvariant subspaces
for some classes of operators will be given in sections that follow.

\begin{proposition} \label{prop13}
Let $A\in B(\eX)$ and let $U\in B(\eX)$ be an invertible operator. A closed subspace $\eM\subseteq \eX$ is an ultrainvariant subspace of $A$
if and only if $U\eM$ is an ultrainvariant subspace of $UAU^{-1}$.
\end{proposition}

\begin{proof}
Assume that $\eM$ is ultrainvariant. By \eqref{eq06}, every operator in $\LCom(UAU^{-1};U\eM)$
is of the form $UTU^{-1}$ for some $T\in \LCom(A;\eM)$. Hence, for every $x\in \eM$, we have $(UTU^{-1})(Ux)=UTx\in U\eM$, that is, $U\eM$ is invariant
for every operator from $\LCom(UAU^{-1};U\eM)$ and is therefore ultrainvariant. It is clear that the opposite implication holds, as well.
\end{proof}

Recall that a closed subspace $\eX_1\subseteq \eX$ is reducing for $A\in B(\eX)$ if $\eX_1\in \Lat(A)$ and there exists $\eX_2\in \Lat(A)$ such that $
\eX=\eX_1\oplus \eX_2$. In this case we will say that $(\eX_1,\eX_2)$ is a reducing pair for $A$. Note that $A=A_1\oplus A_2$ with respect to the 
decomposition $\eX=\eX_1\oplus \eX_2$. Let $P_j:\eX \to \eX_j$ $(j=1,2)$ be given by $P_1(x_1\oplus x_2)=x_1$
and $P_2(x_1\oplus x_2)=x_2$.

\begin{proposition} \label{prop16}
Let $A\in B(\eX)$ and let $(\eX_1,\eX_2)$ be a reducing pair for $A$. If $\eM$ is an ultrainvariant subspace of $A$, then $\eM_j=P_j\eM$ is an ultrainvariant 
subspace of $A_j$ $(j=1,2)$.
\end{proposition}

\begin{proof}
Assume that $\eM=\eM_1\oplus \eM_2$ is ultrainvariant. Let $T_1\in \LCom(A_1;\eM_1)$ be arbitrary and let $T\in B(\eX)$ be such that $T=T_1 \oplus 0$ with respect to the
decomposition $\eX=\eX_1\oplus \eX_2$. For an arbitrary $x=x_1\oplus x_2\in \eM$, we have $TAx=T_1A_1x_1\oplus 0=A_1T_1x_1\oplus 0=ATx$ which means that $T\in \LCom(A;\eM)$.
Since $\eM$ is ultrainvariant it is invariant for $T$. In particular, for every $x_1\in \eM_1$ we have $T(x_1\oplus 0)=T_1x_1\oplus 0\in \eM_1\oplus \eM_2$,
that is, $\eM_1$ is invariant for $T_1\in \LCom(A_1;\eM_1)$. We have seen that $\eM_1$ is an ultrainvariant subspace of $A_1$ and a similar proof
shows that $\eM_2$ is an ultrainvariant subspace of $A_2$.
\end{proof}

Let $n\geq 2$ be an integer. For every $j\in \{ 1, \ldots, n\}$, let $\eX_j$ be a Banach space and $A_j\in B(\eX_j)$. Denote $\eX=\eX_1\oplus \cdots \oplus \eX_n$
and $A=A_1\oplus \cdots\oplus A_n$. Let $\eM_j\subseteq \eX_j$ be a complemented closed subspace and let $\eN_j\subseteq \eX_j$ be a closed subspace such that
$\eX_j=\eM_j \oplus \eN_j$. Let $\left[ \begin{smallmatrix} A_{11}^{(j)} & A_{12}^{(j)}\\ A_{21}^{(j)} & A_{22}^{(j)}\end{smallmatrix}\right]$ be the operator matrix
of $A_j$ with respect to the decomposition $\eX_j=\eM_j\oplus \eN_j$.

\begin{proposition} \label{prop17}
If, for every $j\in \{ 1, \ldots, n\}$, $\eM_j$ is an ultrainvariant subspace of $A_j$ and $\sigma(A_{11}^{(i)})\cap \sigma(A_j)=\emptyset$ whenever $i\ne j$, then
$\eM:=\eM_1\oplus \cdots \oplus \eM_n$ is an ultrainvariant subspace of $A$.
\end{proposition}

\begin{proof}
It is enough to consider the case $n=2$. We will identify $\eX$ with $\eM_1\oplus \eN_1\oplus \eM_2\oplus \eN_2$.
 Then $\eM=\eM_1\oplus \{ 0\}\oplus \eM_2\oplus \{ 0\}$. With respect to the
decomposition $\eX=\eM_1\oplus \eN_1\oplus \eM_2\oplus \eN_2$, operator $A$ has operator matrix
$$ \left[ \begin{smallmatrix} 
A_{11}^{(1)} & A_{12}^{(1)} & 0 & 0\\
0 & A_{22}^{(1)} & 0 & 0\\
0 & 0 & A_{11}^{(2)} & A_{12}^{(2)}\\
0 & 0 & 0 & A_{22}^{(2)}
\end{smallmatrix}
 \right]. $$
Let $T\in \LCom(A;\eM)$ be arbitrary and let $\left[ T_{ij}\right]_{i,j=1}^{4}$ be its operator matrix with respect to the decomposition
$\eX=\eM_1\oplus \eN_1\oplus \eM_2\oplus \eN_2$. Hence, $TAx=ATx$ for every $x=x_1\oplus 0\oplus x_2\oplus 0\in \eM$.
If $x=x_1\oplus 0 \oplus 0\oplus 0$, where $x_1\in \eM_1$ is arbitrary, then it follows from $TAx=ATx$ that
\begin{equation} \label{eq07}
\begin{split}
T_{11}A_{11}^{(1)}&=A_{11}^{(1)}T_{11}+A_{12}^{(1)}T_{21}\\
T_{21}A_{11}^{(1)}&=A_{22}^{(1)}T_{21}\\
T_{31}A_{11}^{(1)}&=A_{11}^{(2)}T_{31}+A_{12}^{(2)}T_{41}\\
T_{41}A_{11}^{(1)}&=A_{22}^{(2)}T_{41}.
\end{split}
\end{equation}
Let $W_{11}\in B(\eX_1)$ be the operator whose operator matrix with respect to the decomposition $\eX_1=\eM_1\oplus \eN_1$
is $\left[\begin{smallmatrix} T_{11} & T_{12}\\ T_{21} & T_{22}\end{smallmatrix}\right]$. The first two equalities in \eqref{eq07} give
$W_{11}\in \LCom(A_1;\eM_1)$. Since $\eM_1$ is an ultrainvariant subspace of $A_1$ it is invariant for $W_{11}$ and therefore
$T_{21}=0$. Hence, by the first equality in \eqref{eq07}, $T_{11}\in \bigl(A_{11}^{(1)}\bigr)'$. Denote by $W_{21}$ the operator in $B(\eX_1,\eX_2)$
whose operator matrix with respect to the decompositions $\eX_1=\eM_1\oplus \eN_1$ and $\eX_2=\eM_2\oplus \eN_2$ is 
$\left[\begin{smallmatrix} T_{31} & T_{32}\\ T_{41} & T_{42}\end{smallmatrix}\right]$. It follows from the last two equalities in \eqref{eq07}
that column $\left[\begin{smallmatrix} T_{31} \\ T_{41}\end{smallmatrix}\right]$ intertwines $A_{11}^{(1)}$ and $A_2$. Since, by the assumption,
$\sigma(A_{11}^{(1)})\cap \sigma(A_2)=\emptyset$ we conclude that $T_{31}=0$ and $T_{41}=0$.

Let now $x_2\in \eM_2$ be arbitrary and $x=0\oplus 0\oplus x_2 \oplus 0$. A similar reasoning as in the previous paragraph shows
that $TAx=ATx$ if and only if $T_{13}=0$, $T_{23}=0$, $T_{43}=0$ and $T_{33}\in \bigl(A_{11}^{(2)}\bigr)'$. Hence, the operator matrix of $T$
with respect to the decomposition $\eX=\eM_1\oplus \eN_1\oplus \eM_2\oplus \eN_2$ is
$$ \left[ \begin{smallmatrix}
T_{11} & T_{12} & 0 & T_{14}\\
0 & T_{22} & 0 & T_{24}\\
0 & T_{32} & T_{33} & T_{34}\\
0 & T_{42} & 0 & T_{44}
\end{smallmatrix}
 \right]. $$
It is not hard to check now that $\eM=\eM_1\oplus \{ 0\}\oplus \eM_2\oplus \{ 0\}$ is invariant for $T$. Since $T\in \LCom(A;\eM)$ was arbitrary 
we conclude that $\eM$ is an ultrainvariant subspace of $A$.
\end{proof}

\begin{proposition} \label{prop05}
Let $A\in B(\eX)$ and let $\eM\subseteq \eX$ be a closed subspace. Then 
$$\LCom\bigl(A;\LCom(A;\eM)\eM\bigr)\LCom(A;\eM)\eM$$
is the smallest ultrainvariant subspace of $A$ which contains $\eM$.
\end{proposition}

\begin{proof}
By Corollary \ref{cor01}, $\LCom\bigl(A;\LCom(A;\eM)\eM\bigr)$ is the largest algebra which is contained in $\LCom(A;\eM)$. It follows, by Theorem \ref{theo04},
that $\bigl(\LCom(A;\eM)\eM\bigr)_A=\LCom\bigl(A;\LCom(A;\eM)\eM\bigr)\LCom(A;\eM)\eM$ and
$\bigl(\LCom(A;\eM)\eM\bigr)_A$ is an ultrainvariant subspace of $A$. It is clear that $\eM\subseteq \bigl(\LCom(A;\eM)\eM\bigr)_A$.

Suppose that $\eK\subseteq \eX$ is an ultrainvariant subspace of $A$ such that $\eM \subseteq \eK$. Then
$\LCom(A;\eK)$ is an algebra and $\LCom(A;\eM)\supseteq \LCom(A;\eK)$.
Since $\LCom\bigl(A;\bigl(\LCom(A;\eM)\eM\bigr)_A\bigr)=\LCom\bigl(A;\LCom(A;\eM)\eM\bigr)$
is the largest algebra contained in $\LCom(A;\eM)$ we have $\LCom(A;\eK)\subseteq \LCom\bigl(A;\LCom(A;\eM)\eM\bigr)$. Hence, if $T\in \LCom(A;\eK)$, then it commutes
with $A$ at every vector from $\LCom(A;\eM)\eM$ which means that $\LCom(A;\eM)\eM\subseteq \eK_{_A}=\eK$ and therefore
$\bigl(\LCom(A;\eM)\eM\bigr)_{_A}\subseteq \eK$.
\end{proof}

For $A\in B(\eX)$, let $\Latu(A)$ denote the family of all ultrainvariant subspaces of $A$.

\begin{proposition} \label{prop12}
Let $A\in B(\eX)$ and let $\{\eM_j;\; j\in J\}$ be an arbitrary family of ultrainvariant subspaces of $A$. Then $\bigvee\limits_{j\in J}\eM_j$
and $\bigcap\limits_{j\in J}\eM_j$ are ultrainvariant subspaces of $A$. Hence, $\Latu(A)$ is a sublattice of $\Lath(A)$.
\end{proposition}

\begin{proof}
Denote $\eK=\bigvee\limits_{j\in J}\eM_j$. 
By \eqref{eq05}, we have $\LCom(A;\eK)=\bigcap\limits_{j\in J}\LCom(A;\eM_j)$. Hence, $\eM_j$ $(j\in J)$ is invariant for an arbitrary 
$T\in \LCom(A;\eK)$. Assume that $x\in \eK$ is such that $x=x_1+\cdots+x_k$, where $x_i\in \eM_{j_i}$ $(i=1, \ldots, k)$. Then $Tx=Tx_1+\cdots+Tx_k
\in \eM_{j_1}+\cdots+\eM_{j_k}\subseteq \eM$. It follows that $T\eK\subseteq \eK$. Since $T\in \LCom(A;\eK)$ is arbitrary we conclude that $\eK$
is ultrainvariant.

Denote $\eL=\bigcap\limits_{j\in J}\eM_j$ and let $\eL'\subseteq \eX$ be the smallest ultrainvariant subspace of $A$ such that $\eL \subseteq \eL'$ (see Proposition \ref{prop05}). Since, for every $j\in J$,
$\eL \subseteq \eM_j$ and $\eM_j$ is ultrainvariant we have $\eL'\subseteq \eM_j$. Hence $\eL'\subseteq \bigcap\limits_{j\in J}\eM_j=\eL$, that is, $\eL=\eL'$, i.e., $\eL$ is an ultrainvariant subspace of $A$.
\end{proof}

The following example shows that even when $\eM \in \Lat(A)$ is the girder of $\LCom(A;\eM)$ it  is not necessary an
ultrainvariant subspace of $A$, that is, $\LCom(A;\eM)$ is not necessary an algebra.

\begin{example} \label{ex02}
Let $\eX=\eX_1\oplus \eX_2\oplus \eX_3$, where $\eX_j$ ($j=1,2,3$) are non-trivial subspaces of $\eX$. Let $A$ be the projection on $\eX_1\oplus \eX_2$
along $\eX_3$ and let $\eM=\eX_1\oplus \{ 0\} \oplus \eX_3$, that is, in $\eM$ are vectors of the form $x_1\oplus 0\oplus x_3$, where $x_1\in \eX_1$ and $x_3\in \eX_3$
are arbitrary. Hence, $\eM$ is a non-trivial subspace of $\eX$ and $\eM\in \Lat(A)$. 

Let $T\in \LCom(A;\eM)$ be arbitrary and let $\left[ T_{ij}\right]_{i,j=1}^{3}$ be its operator matrix with respect to the decomposition $\eX=\eX_1\oplus \eX_2\oplus \eX_3$.
For every $x=x_1\oplus 0\oplus x_3\in \eM$, we have $TAx=T_{11}x_1\oplus T_{21}x_1\oplus T_{31}x_1$ and $ATx=(T_{11}x_1+T_{13}x_3)\oplus (T_{21}x_1+T_{23}x_3)\oplus 0$.
It follows from $TAx=ATx$ that $T_{11}x_1=T_{11}x_1+T_{13}x_3$, $T_{21}x_1=T_{21}x_1+T_{23}x_3$ and $T_{31}x_1=0$ for all $x_1\in \eX_1$ and $x_3\in \eX_3$.
Thus, $T_{13}=0$, $T_{23}=0$ and $T_{31}=0$. An operator $T$ is in $\LCom(A;\eM)$ if and only if its operator matrix with respect to the decomposition
$\eX=\eX_1\oplus \eX_2\oplus \eX_3$ is of the form $\left[ \begin{smallmatrix} T_{11} & T_{12} & 0\\ T_{21} & T_{22} & 0\\ 0 & T_{32} & T_{33}\end{smallmatrix}\right]$,
where $T_{ij}\in B(\eX_j, \eX_i)$ are arbitrary. It follows that $\LCom(A;\eM)$ is not an algebra.

To determine $\eM_A$, let $x=x_1\oplus x_2\oplus x_3\in \eX$ be such that $TAx=ATx$ for every $T\in \LCom(A;\eM)$. Since $TAx=(T_{11}x_1+T_{12}x_2)\oplus(T_{21}x_1+T_{22}x_2)\oplus T_{32}x_2$
and $ATx=(T_{11}x_1+T_{12}x_2)\oplus(T_{21}x_1+T_{22}x_2)\oplus 0$ we see that $T_{32}x_2=0$ for every $T_{32}\in B(\eX_2,\eX_3)$. As the involved subspaces are
non-trivial we conclude that $x_2=0$. This proves that $\eM_A=\eM$. However, $\eM$ is not an ultrainvariant subspace of $A$ because $\LCom(A;\eM)$ is not an algebra.
Moreover, $\LCom(A;\eM)\eM=\eX$. Namely, let $x=x_1\oplus x_2\oplus x_3$ be an arbitrary vector from $\eX$. Let $y=y_1\oplus 0\oplus y_3\in \eM$ be such that
$y_1\ne 0$ and $y_3\ne 0$. Then there exist $T_{11}\in B(\eX_1)$, $T_{21}\in B(\eX_1,\eX_2)$ and $T_{33}\in B(\eX_3)$ such that $T_{11}y_1=x_1$,
$T_{21}y_1=x_2$ and $T_{33}y_3=x_3$. It is clear that $T$ whose operator matrix with respect to the decomposition $\eX=\eX_1\oplus \eX_2\oplus \eX_3$
is $\left[ \begin{smallmatrix} T_{11} & 0 & 0\\ T_{21} & 0 & 0\\ 0 & 0 & T_{33}\end{smallmatrix}\right]$ is in $\LCom(A;\eM)$ and $Ty=x$.~\hfill $\Box$
\end{example}

%************************************************************************
%************************************************************************
\section{Operators with non-trivial ultrainvariant subspaces} \label{sec05}
\setcounter{theorem}{0}
%************************************************************************
%************************************************************************

Let $A\in B(\eX)$ and assume that $\eM\in \Lat(A)$ is complemented, say $\eX=\eM \oplus \eN$. 
Let $\left[\begin{smallmatrix} A_{11} & A_{12}\\ 0 & A_{22}\end{smallmatrix}\right]$ be the operator matrix of $A$ with respect this decomposition. Next proposition is basically a reformulation of Corollary \ref{cor05}.

\begin{proposition} \label{prop18}
If $\sigma(A_{11})\cap \sigma(A_{22})=\emptyset$, then $\eM$ is an ultrainvariant subspace of $A$.
\end{proposition}

\begin{proof}
If $\eM$ is trivial, then there is nothing to prove. Assume therefore that $\{ 0\} \ne \eM \ne \eX$.
 By Corollary \ref{cor05},  $\LCom(A;\eM)$ is an algebra. Hence, if $\left[ \begin{smallmatrix} T_{11} & T_{12}\\ T_{21} & T_{22}\end{smallmatrix} \right]$ is the operator matrix of $T\in \LCom(A;\eM)$ with respect to the decomposition $\eX=\eM\oplus \eN$, then $T_{11}\in (A_{11})'$ and $T_{21}=0$ by Theorem \ref{theo07}. It follows that $\eM$ is invariant for every $T\in \LCom(A;\eM)$.
\end{proof}

Recall that the ascent $\alpha(A)$ of $A\in B(\eX)$ is the smallest positive integer $k$ such that $\Ker(A^k)=\ker(A^{k+1})$. If there is no positive
integer with this property, then $\alpha(A)=\infty$. The descent $\delta(A)$ of $A$ is the smallest positive integer $k$ such that $\Image(A^k)=\Image(A^{k+1})$
and $\delta(A)=\infty$ if there is no $k$ with this property. Assume that $A$ has finite ascent and descent. Then $\alpha(A)=n=\delta(A)$ for a positive
integer $n$ and, by \cite[Lemma 2.21]{AA}, $\Image(A^n)$ is a closed complement of $\Ker(A^n)$, that is $\eX=\Ker(A^n)\oplus \Image(A^n)$. Hence,
$(\Ker(A^n),\Image(A^n))$ is a reducing pair of $A$ and therefore $A=A_1\oplus A_2$, where $A_1\in B(\Ker(A^n))$ and $A_2\in B(\Image(A^n))$. By \cite[Theorem 2.23]{AA},
$A_1$ is nilpotent and $A_2$ is invertible which implies $\sigma(A_1)\cap \sigma(A_2)=\emptyset$. It is clear now that the following corollary is a simple
consequence of Proposition \ref{prop18}

\begin{corollary} \label{cor09}
If $A\in B(\eX)$ has finite ascent and descent $n$, then $\Ker(A^n)$ and $\Image(A^n)$ are ultrainvariant subspaces of $A$.
\end{corollary}

It is not hard to see that for every $\lambda \in \bC$ the kernel $\Ker(A-\lambda \IX)$ is an ultrainvariant subspace of $A\in B(\eX)$. Namely, 
if $T\in \LCom(A;\Ker(A-\lambda \IX))$, then $ATx=TAx=\lambda Tx$ for every $x\in \Ker(A-\lambda \IX)$, that is, $\Ker(A-\lambda \IX)$ is invariant for every 
$T\in \LCom(A;\Ker(A-\lambda \IX))$. Hence, if $A$ is not a scalar multiple of $\IX$ and the point spectrum $\sigma_p(A)$ is not empty, then $A$ has a 
non-trivial ultrainvariant subspace. More can be said. We need the following lemma.

\begin{lemma} \label{lem03} 
Let $A\in B(\eX)$ and let $\pA \subseteq B(\eX)$ be the strongly closed subalgebra generated by $A$ and $\IX$. If $\eM \in \Lat(A)$, then 
$\LCom(A;\eM)\subseteq \LCom(B;\eM)$ for every $B\in \pA$.
\end{lemma}

\begin{proof}
Assume first that $B=A-\lambda \IX$, where $\lambda\in \bC$ is such that $B$ is an invertible operator. Hence, $B\eM=\eM$.
It is clear that $C(A;\eM)=\LCom(B;\eM)$. Let us check that $\LCom(B;\eM)\subseteq \LCom(B^k;\eM)$ for every integer $k\geq 0$. 
If $k=0$ or $k=1$, then the assertion is trivial. Since $B$ is invertible and $B\eM=\eM$ we have $\LCom(B;\eM)=\LCom(B;\eM)B$, by Proposition \ref{prop04}. 
Hence, if $T\in \LCom(B;\eM)$, then $TB\in \LCom(B;\eM)$ and, by induction, $TB^j\in \LCom(B;\eM)$ for every $j\geq 2$. Let $k\geq 2$ and assume that
$TB^jx=B^jTx$ holds for every $0\leq j \leq k-1$ and every $x\in \eM$. Then
$TB^kx=B^{k-1}TBx$ and $TBx=BTx$ for every $x\in \eM$. It follows that $TB^kx=B^kTx$ for every $x\in \eM$. Since $A=B+\lambda \IX$ we conclude that
$$TA^kx=T\sum_{j=0}^{k}{k\choose j}\lambda^j B^{k-j}x=\sum_{j=0}^{k}{k\choose j}\lambda^j B^{k-j}Tx=A^kTx\qquad \text{for every}\; x\in \eM.$$
It follows that $Tp(A)x=p(A)Tx$ for every $x\in \eM$ and every polynomial $p\in \bC[z]$. Let $B\in \pA$ be arbitrary. Let $T\in \LCom(A;\eM)$ and $x\in \eM$.
For every $\varepsilon>0$, there exists a polynomial $p$ such that $\| Bx-p(A)x\|< \varepsilon$ and $\| BTx-p(A)Tx\|< \varepsilon$. 
Hence $\| TBx-BTx\|\leq \| T\| \| Bx-p(A)x\|+ \| BTx-p(A)Tx\|<\varepsilon (\| T\|+1)$, which gives $TBx=BTx$.
\end{proof}

\begin{proposition} \label{prop03}
Let $A\in B(\eX)$ and let $\pA \subseteq B(\eX)$ be the strongly closed subalgebra generated by $A$ and $\IX$. If $B\in \pA$ and $\eM\in \Lat(A)$ is an
ultrainvariant subspace of $B$, then $\eM$ is an ultrainvariant subspace of $A$. In particular, every kernel $\Ker(B-\lambda \IX)$ $(\lambda \in \bC)$
is an ultrainvariant subspace of $A$.
\end{proposition}

\begin{proof}
By Lemma \ref{lem03}, $\LCom(A;\eM)\eM\subseteq \LCom(B;\eM)\eM=\eM$. Thus, $\eM$ is invariant for every operator from $\LCom(A;\eM)$, that is,
it is an ultrainvariant subspace of $A$.
\end{proof}

It follows from the following proposition that a reducing subspace does not need to be ultrainvariant.

\begin{proposition} \label{prop07}
Let $A\in B(\eX)$. Assume that $(\eM, \eN)$ is a pair of reducing subspaces for $A$. Let $A$ be equal to $A_{_{\eM}}\oplus A_{_{\eN}}$ with respect to the 
decomposition $\eX=\eM \oplus \eN$. Then $\eM$ is an ultrainvariant subspace of $A$ if and only if $\LInt(A_{_{\eM}},A_{_{\eN}})=\{ 0\}$.
\end{proposition}

\begin{proof}
Suppose that $\eM$ is ultrainvariant and let $S\in \LInt(A_{_\eM},A_{_\eN})$ be arbitrary. Then $T(x\oplus y)=0\oplus Sx$, where $x\in \eM, y\in \eN$
are arbitrary, defines a bounded linear operator on $\eX$. Since
$AT(x\oplus 0)=0\oplus A_{_\eN}Sx=0\oplus SA_{_\eM}x=TA(x\oplus 0)$ for all $x\in \eM$
we have $T\in \LCom(A;\eM)$. It follows from $T\eM\subseteq \eM$ that $Sx=0$ for every $x\in \eM$, that is, $S=0$.

Assume now that $\LInt(A_{_\eM},A_{_\eN})=\{ 0\}$. Let $T\in \LCom(A;\eM)$ be arbitrary and let $\left[\begin{smallmatrix} T_{11} & T_{12}\\ T_{21} & T_{22}\end{smallmatrix} \right]$
be its operator matrix with respect to the decomposition $\eX=\eM\oplus \eN$. It follows from $TA(x\oplus 0)=AT(x\oplus 0)$, where $x\in \eM$ is arbitrary,
that $T_{21}\in \LInt(A_{_\eM},A_{_\eN})$, that is, $T_{21}=0$. Hence, $\eM$ is an ultrainvariant subspace of $A$.
\end{proof}

Let $A\in B(\eX)$. A subset $\emptyset \subsetneq \sigma \subsetneq \sigma(A)$ is an isolated part of $\sigma(A)$ if both $\sigma$ and $\sigma(A)\setminus \sigma$
are closed sets. It is well-known, see \cite[\S I.2]{GGK1}, that for an isolated part $\sigma$ of $\sigma(A)$ there exists an idempotent $P_\sigma\in B(\eX)$, 
the Riesz projection of $A$ corresponding to $\sigma$, such that $\eM=\Image(P_\sigma)$ and $\eN=\Image(I-P_\sigma)$ are in $\Lat(A)$ and 
$\sigma(A|_\eM)=\sigma$, $\sigma(A|_\eN)=\sigma(A)\setminus \sigma$. Hence, the following result is an immediate consequence of Proposition \ref{prop07}.

\begin{corollary} \label{cor06}
Let $A\in B(\eX)$. If $\sigma$ is an isolated part of $\sigma(A)$, then $\Image(P_\sigma)$ is an ultrainvariant subspace of $A$.
\end{corollary}

Corollary \ref{cor06} is a special case of the next result for which we need some notions from the local spectral theory.
Let $A\in B(\eX)$ and $x\in \eX$. The local resolvent set $\rho_A(x)$ of $A$ at $x$ is the union of all open subsets $U\subseteq \bC$ for which there
exists an analytic function $f:U\to \eX$ such that $(A-\lambda \IX)f(\lambda)=x$ for all $\lambda\in U$. The local spectrum of $A$ at $x$ is then defined
as $\sigma_A(x)=\bC\setminus \rho_A(x)$. It is obvious that $\sigma_A(x)$ is a closed subset of the spectrum $\sigma(A)$. For an arbitrary $F\subseteq \bC$,
the local spectral subspace of $A$ corresponding to $F$ is $\eX_A(F)=\{ x\in \eX;\; \sigma_A(x)\subseteq F\}$. By \cite[Proposition 1.2.16 (a)]{LN},
$\overline{\eX_A(F)}$ is a hyperinvariant subspace of $A$.

\begin{proposition} \label{prop02}
$\overline{\eX_A(F)}$ is an ultrainvariant subspace of $A$.
\end{proposition}

\begin{proof}
Let $T\in \LCom(A;\overline{\eX_A(F)})$ be arbitrary. Let $x\in \eX_A(F)$ and let $\lambda_0\in \rho_A(x)$. Then there is an open neighborhood $U\subseteq 
\rho_A(x)$ of $\lambda_0$ and an analytic function $f:U\to \eX$ such that $(A-\lambda \IX)f(\lambda)=x$ for all $\lambda\in U$. By \cite[Lemma 1.2.14 ]{LN},
$\sigma_A(f(\lambda))=\sigma_A(x)$ for every $\lambda \in U$. Hence, $f(\lambda)\in \eX_A(F)$ for all $\lambda \in U$. It follows that
$(A-\lambda \IX)Tf(\lambda)=T(A-\lambda \IX)f(\lambda)=Tx$ for all $\lambda \in U$. Since $\lambda \mapsto Tf(\lambda)$ is an analytic function on $U$
we conclude that $\lambda_0\in \sigma_A(Tx)$ and therefore $\sigma_A(Tx)\subseteq \sigma_A(x)\subseteq F$. Thus, $\overline{\eX_A(F)}$ is invariant for $T$.
\end{proof}

Recall that $A\in B(\eX)$ is a decomposable operator if for every open cover $\bC=U\cup V$ there exist spaces $\eM, \eN\in \Lat(A)$ such that $\eX=\eM+\eN$
and $\sigma(A|_{\eM})\subseteq U$, $\sigma(A|_{\eN})\subseteq V$ (see Definition 1.1.1 in \cite{LN}). If $A$ is decomposable, then the local spectral subspace
$\eX_A(F)$ is closed whenever $F$ is a closed subset of $\bC$ (see Definition 1.2.18 and Theorem 1.2.29 in \cite{LN})
and $\sigma(A)=\cup\{ \sigma_A(x);\; x\in \eX\}$ (see Proposition 1.3.2 in \cite{LN}).

\begin{corollary} \label{cor02}
If $A\in B(\eX)$ is a decomposable operator such that the spectrum $\sigma(A)$ contains at least two points, then $A$ has a proper non-trivial
ultrainvariant subspace. 
\end{corollary}

\begin{proof}
Let $\lambda_1\ne \lambda_2$ be points in $\sigma(A)$ and let $\bC=U\cup V$ be an open cover such that $\lambda_1\in U$ and $\lambda_2\not\in \overline{U}$.
Then $\eX_A(\overline{U})$ and $\eX_A(\overline{V})$ are non-trivial and proper closed subspaces of $\eX$. Indeed, by \cite[Theorem 1.2.23]{LN}, 
$\eX=\eX_A(\overline{U})+\eX_A(\overline{V})$.
Since $\sigma_A(x)\subseteq \overline{U}\cap \sigma(A)\subsetneq \sigma(A)$ and $\sigma(A)=\cup\{ \sigma_A(x);\; x\in \eX\}$ there are vectors which are
not in $\eX_A(\overline{U})$. Similarly, there are vectors which are not in $\eX_A(\overline{V})$. Hence, 
 $\eX_A(\overline{U})$ and $\eX_A(\overline{V})$ are proper and non-trivial. 
By Proposition \ref{prop02}, these subspaces are ultrainvariant. 
\end{proof}

Assume that $A\in B(\eX)$ is a nilpotent operator of order $n\geq 2$, that is, $A^n=0$ and $A^{n-1}\ne 0$. 
It is well-known that
\begin{equation*}
\{ 0\}=\Ker(A^0)\subseteq \Ker(A)\subseteq \Ker(A^2)\subseteq \cdots \subseteq\Ker(A^{n-1})\subseteq \Ker(A^n)=\eX
\end{equation*}
and
\begin{equation*}
\{ 0\}=\Image(A^n)\subseteq \overline{\Image(A^{n-1})}\subseteq \overline{\Image(A^{n-2})}\subseteq \cdots \subseteq\overline{\Image(A)}\subseteq \Image(A^0)=\eX
\end{equation*}
are two chains of hyperinvariant subspaces of $A$. 
By \cite[Remarque p. 317]{Bar}, subspaces $\overline{\Image(A^{n-1})}$ and
$\Ker(N^{n-1})$ are the smallest, respectively the largest, non-trivial hyperinvariant subspaces of $A$.
For every $j=0,1,\ldots, n-1$, we have $\overline{\Image(A^{n-j})}\subsetneq \Image(A^{n-j-1})$ and $\Ker(A^j)\subsetneq \Ker(A^{j+1})$.
Of course, $\overline{\Image(A^{n-j})}\subseteq \Ker(A^j)$. Later we will need the following simple facts.

\begin{lemma} \label{lem04}
(i) For every $j=1, \ldots, n-1$, there exists $u\in \Image(A^{n-j})$ such that $u\not\in \Ker(A^{j-1})$.

(ii) Let $e\in \eX$ be such that $A^k e\ne 0$ for an integer $k\in \{0,1, \ldots, n-1\}$. Then
$A^k e\not\in \bigvee\{ e, Ae,\ldots, A^{k-1}e\}$, that is, $e, Ae,\ldots, A^k e$ are linearly independent.
\end{lemma}

\begin{proof}
(i) If $\Image(A^{n-j})\subseteq \Ker(A^{j-1})$, then we would have $A^{j-1}(N^{n-j}x)=0$ for every $x\in \eX$, that is, it would be $A^{n-1}=0$.

(ii) Towards a contradiction assume that there exist numbers $\alpha_0, \ldots, \alpha_{k-1}$ such that
$$ A^ke=\alpha_0 e+\alpha_1 Ae+\cdots+\alpha_{k-1} A^{k-1}e. $$
Since $A^k e\ne 0$ there exists the smallest index $0\leq j\leq k-1$ such that $\alpha_j\ne 0$. Hence,
$$ A^ke=\alpha_j A^j e+\alpha_{j+1} A^{j+1}e+\cdots+\alpha_{k-1} A^{k-1}e. $$
Let $m$ be the integer such that $A^me\ne 0$ and $A^{m+1}e=0$.
Multiply the above equality by $A^{m-j}$. Since $j<k$ we have
$$ 0=A^{m-j+k}e=\alpha_jA^me+\alpha_{j+1}A^{m+1}+\cdots+\alpha_{k-1}A^{m-j+k-1}e=\alpha_j A^m e.$$
It follows that $\alpha_j=0$ as $A^m e\ne 0$. This is a contradiction.
\end{proof}

\begin{lemma} \label{lem01}
Let $j\in \{1, \ldots, n\}$. If $e\in \Ker(A^j)$ and $\xi\in \eX^*$, then
$$(e\otimes \xi) A^{j-1}+A(e\otimes \xi) A^{j-2}+\cdots+A^{j-1}(e\otimes \xi) \in \LCom(A;\overline{\Image(A^{n-j})}).$$
\end{lemma}

\begin{proof}
It is not hard to check that
$(e\otimes \xi) A^{j-1}+A(e\otimes \xi) A^{j-2}+\cdots+A^{j-1}(e\otimes \xi)$ and $A$ commute locally at
every vector $A^{n-j}x\in\Image(A^{n-j})$, where $x\in \eX$ is arbitrary.
\end{proof}

\begin{theorem} \label{theo03}
Let $A\in B(\eX)$ be a nilpotent operator of order $n\geq 2$. A closed subspace $\eM \subseteq \eX$ is ultrainvariant for $A$ if and only if
$\eM=\Ker(A^j)$ for some $j\in \{ 0,1,\ldots,n\}$.
\end{theorem}

\begin{proof}
By Proposition \ref{prop03}, every kernel $\Ker(A^j)$ $(j= 0,1,\ldots,n)$ is an ultrainvariant subspace of $A$. For the opposite implication,
suppose that $\eM$ is an ultrainvariant subspace for $A$. Then $\eM$ is hyperinvariant for $A$ and therefore, by \cite[Lemma 5]{Bar}, there
exists a unique $j\in \{ 0,1, \ldots, n\}$ such that $\overline{\Image(A^{n-j})}\subseteq \eM \subseteq \Ker(A^j)$.

We claim that $\LCom(A,\overline{\Image(A^{n-j})})\,\overline{\Image(A^{n-j})}=\Ker(A^j)$, for every $j\in \{0, 1, \ldots, n\}$.  
If $j=0$ or $j=n$, then the assertion is trivial. Assume therefore that $1\leq j \leq n-1$. Let $T\in \LCom(A,\overline{\Image(A^{n-j})})$
and $x\in \Image(A^{n-j})$ be arbitrary. Let $z\in \eX$ be such that $x=A^{n-j}z$. By Lemma \ref{lem03}, $A^j Tx=TA^j x=TA^nz=0$.
This proves that $\LCom(A,\overline{\Image(A^{n-j})})\,\overline{\Image(A^{n-j})}\subseteq \Ker(A^j).$
For the opposite inclusion, assume that $x\in \Ker(A^j)$. Let $e\in \Image(A^{n-j})$ be such that $A^{j-1}e\ne 0$ (it exists by Lemma \ref{lem04}~(i)).
Since, by Lemma \ref{lem04}~(ii), vectors $e, Ae, \ldots, A^{j-1}e$ are linearly independent there exists $\xi \in \eX^*$ such that
$\langle A^{j-1}e,\xi\rangle=1$ and $\langle A^{i}e,\xi\rangle=0$ for $i=0,\ldots, j-2$. Denote $T=(x\otimes \xi)A^{j-1}+A(x\otimes \xi)A^{j-2}+\cdots+A^{j-1}(x\otimes \xi)$.
By Lemma \ref{lem01}, $T\in \LCom(A,\overline{\Image(A^{n-j})})$. Since $Te=x$ the equality $\LCom(A,\overline{\Image(A^{n-j})})\,\overline{\Image(A^{n-j})}=\Ker(A^j)$
is proved. By Proposition \ref{prop05}, 
$$\LCom\bigl(A;\LCom(A,\overline{\Image(A^{n-j})})\overline{\Image(A^{n-j})}\bigr)\,\LCom(A,\overline{\Image(A^{n-j})})\overline{\Image(A^{n-j})}
=\LCom\bigl(A;\Ker(A^j)\bigr)\,\Ker(A^j)=\Ker(A^j)$$
is the smallest ultrainvariant subspace of $A$ which contains $\overline{\Image(A^{n-j})}$. Hence, $\eM=\Ker(A^j)$.
\end{proof}

An operator $A\in B(\eX)$ is  algebraic if there exists a non-zero polynomial $p$ such that $p(A)=0.$ It is not hard to see that for an algebraic operator 
$A$ there exists a unique monic polynomial $q_A$ of the minimal degree, called the  minimal polynomial of $A$, such that $q_A(A)=0$. Every operator on a
finite dimensional Banach space is algebraic.
Let $ q_A(z)=(z-\lambda_1)^{n_1}\cdots (z-\lambda_k)^{n_k} $ with $\lambda_1, \ldots, \lambda_k$ distinct. For every $j=1, \ldots, k$, let 
$\eX_j=\Ker\bigl((A-\lambda_j \IX)^{n_i}\bigr)$.
Then $\eX=\eX_1 \oplus \cdots \oplus \eX_k$ and with respect to this decomposition $ A=(\lambda_1 I_{\eX_1}+A_1)\oplus \cdots \oplus (\lambda_kI_{\eX_k}+A_k)$,
where $A_j\in B(\eX_j)$ is a nilpotent operator of order $n_j$ (if $n_j=1$, then $A_j=0$).
The reader is referred to \cite[\S 5.9]{Tay} for details about algebraic operators on a complex Banach space.

\begin{corollary} \label{cor08}
Let $A\in B(\eX)$ be algebraic and let $ A=(\lambda_1I_{\eX_1}+A_1)\oplus \cdots \oplus (\lambda_kI_{\eX_k}+A_k)$ be its decomposition as described above.
Then a closed subspace $\eM\subseteq \eX$ is an ultrainvariant subspace of $A$ if and only if $\eM=\Ker(A_{1}^{m_1})\oplus\cdots\oplus \Ker(A_{k}^{m_k})$
for some $0\leq m_j \leq n_j$ $(j=1, \ldots,k)$.
\end{corollary}

\begin{proof}
Assume that $\eM$ is an ultrainvariant subspace of $A$. Let $\eM=\eM_1\oplus\cdots \oplus \eM_k$, where $\eM_j\subseteq \eX_j$ $(j=1, \ldots, k)$. Since every
subspace $\eX_j$ is reducing subspace of $A$, it follows, by Proposition \ref{prop16}, that $\eM_j$ is an ultrainvariant subspace of $\lambda_1I_{\eX_j}+A_j$.
By Theorem \ref{theo03}, there exists $m_j\in \{ 0,1, \ldots, n_j\}$ such that $\eM_j=\Ker(A_{j}^{m_j})$. The opposite implication follows by Proposition
\ref{prop17} because $\sigma(\lambda_iI_{\eX_i}+A_i)\cap \sigma(\lambda_jI_{\eX_j}+A_j)=\emptyset$ whenever $i\ne j$.
\end{proof}

%************************************************************************
%************************************************************************
\section{Ultrainvariant subspaces of operators on a Hilbert space} \label{sec06}
\setcounter{theorem}{0}
%************************************************************************
%************************************************************************

In this section, $\eH$ is the infinite dimensional separable complex Hilbert space. The inner product of vectors $x, y\in \eH$ is denoted
by $(x,y)$. For a closed subspace $\eM \subseteq \eH$, let $\eM^\perp$ be its orthogonal complement. 

Our first result in this section is the complete description of the lattice of ultrainvariant subspaces of a normal operator. It turns out
that $\Latu(A)=\Lath(A)$ whenever $A$ is normal.

Let $A$ be a normal operator on $\eH$ and let $E$ be its spectral measure. For a Borel subset $\sigma\subseteq \sigma(A)$,
operator $E(\sigma)\in B(\eH)$ is an orthogonal projection. Subspace $\eM=\Image(E(\sigma))$ is reducing for $A$ and 
$\eM^\perp=\Image(E(\sigma(A)\setminus \sigma))$ because $E(\sigma(A)\setminus \sigma)=I-E(\sigma)$ (see Section XXXI.7 in \cite{GGK2} for details).

\begin{theorem} \label{theo05}
Let $A\in B(\eH)$ be a normal operator and let $E$ be its spectral measure. A subspace $\eM\subseteq \eH$ is an ultrainvariant subspace of $A$ if and only if
$\eM$ is a hyperinvariant subspace of $A$, that is, if and only if there exists a Borel subset $\sigma\subseteq \sigma(A)$ such that $\eM=\Image(E(\sigma))$.
\end{theorem}

\begin{proof}
By \cite[Theorem 3]{DP}, $\eM\subseteq \eH$ is a hyperinvariant subspace of $A$ if and only if there exists a Borel subset $\sigma\subseteq \sigma(A)$ such that $\eM=\Image(E(\sigma))$.
Hence, it remains to prove that for every Borel subset $\sigma\subseteq \sigma(A)$ the subspace $\eM=\Image(E(\sigma))$ is ultrainvariant.
We can write $A=A_{_{\eM}}\oplus A_{_{\eM^\perp}}$ with respect to the decomposition $\eH=\eM \oplus \eM^\perp$. We have to see that $\LInt(A_{_{\eM}},A_{_{\eM^\perp}})=\{ 0\}$ because then, by Proposition \ref{prop07},  $\eM$ is an ultrainvariant subspace of $A$. There are a few ways how to see that
$\LInt(A_{_{\eM}},A_{_{\eM^\perp}})=\{ 0\}$. We will show it as follows.
Assume that $S\in \LInt(A_{_{\eM}},A_{_{\eM^\perp}})$. Then operator $T\in B(\eH)$ whose operator matrix with respect to the
decomposition $\eH=\eM \oplus \eM^\perp$ is $\left[\begin{smallmatrix} 0 & 0\\ S & 0\end{smallmatrix}\right]$ commutes with $A$. It follows, by 
The Fuglede-Putnam Theorem, that $T$ commutes with $E(\sigma)$ (see Theorems XXXI.7.4 and XXXI.7.9 in \cite{GGK2}, for instance). 
It is easy to see now that $TE(\sigma)=E(\sigma)T$ gives $S=0$.
\end{proof}

Let $1\leq p <\infty$ and let $(X,\mu)$ be a measure space such that $L^p(X,\mu)$ is separable. For $\phi \in L^{\infty}(X,\mu)$, let $M_\phi$ be the 
multiplication by $\phi$ in $L^p(X,\mu)$. By \cite[Theorem 1]{Hua}, a closed subspace $\eM\subseteq L^p(X,\mu)$ is hyperinvariant for $M_\phi$
if and only if $\eM=M_{\chi_{_\sigma}\circ \phi}\bigl(L^p(X,\mu)\bigr)$ for some Borel set $\sigma \subseteq X$. Here $\chi_{_\sigma}$ denotes the characteristic 
function of $\sigma$ and $ M_{\chi_{_\sigma}\circ \phi}$ is the multiplication operator induced by $\chi_{_\sigma}\circ \phi$. Note that $ M_{\chi_{_\sigma}\circ \phi}$ is
an idempotent. Hence $M_{\chi_{_\sigma}\circ \phi}\bigl(L^p(X,\mu)\bigr)$ is complemented and its complement
$M_{\chi_{_{X\setminus\sigma}}\circ \phi}\bigl(L^p(X,\mu)\bigr)$ is a hyperinvariant subspace of $M_\phi$.
A similar reasoning as in the proof of Theorem \ref{theo05} gives that every hyperinvariant subspace of $M_\phi$ is an ultrainvariant subspace.

For a closed subspace $\eM \subseteq \eH$, let $P_{_\eM}$ be the orthogonal projection onto $\eM$. Let $A\in B(\eH)$. A distance on $\Lat(A)$ is
the mapping $\theta$ which is given by $\theta(\eM,\eN)=\| P_{\eM}-P_{\eN}\| $ $(\eM, \eN\in \Lat(A))$. A subspace $\eM \in \Lat(A)$ is said to be
inaccessible if the only continuous mapping $f:[0,1]\to \Lat(A)$ with $f(0)=\eM$ is the constant map $f(t)=\eM$ $(t\in [0,1])$. A simple modification
of the proof of \cite[Theorem 1]{DP} gives the following result.

\begin{theorem} \label{theo06}
Let $A\in B(\eH)$. If $\eM \in \Lat(A)$ is inaccessible, then it is an ultrainvariant subspace of $A$.
\end{theorem}

\begin{proof}
Let $0\ne T\in \LCom(A;\eM)$ be an arbitrary operator. Denote $\Lambda=\{ \lambda \in \bR;\; 0\leq \lambda < \| T\|^{-1}\}$. Then for each $\lambda \in \Lambda$, 
the operator $\IH-\lambda T$ is invertible and therefore $\eM_\lambda=(\IH-\lambda T)\eM$ is a closed subspace of $\eH$. If $y\in \eM_\lambda$,
then there exists $x\in \eM$ such that $y=(\IH-\lambda T)x$. Hence, $Ay=A(\IH-\lambda T)x=(\IH-\lambda T)Ax\in \eM_\lambda$ because $Ax\in \eM$. This shows that
$\eM_\lambda \in \Lat(A)$ for every $\lambda \in \Lambda$. The map $\lambda \mapsto \eM_\lambda$ is continuous (see the proof of Theorem 1 in \cite{DP}).
Since $\eM=\eM_0$ and $\eM$ is inaccessible we have $\eM=\eM_\lambda$ for each $\lambda \in \Lambda$. Thus, $\eM=(\IH-\lambda T)\eM$ for a non-zero $\lambda$
which implies that $\eM$ is invariant for $T$, that is, it is an ultrainvariant subspace of $A$.
\end{proof}

It is clear that $\eM \in \Lat(A)$ is inaccessible if it is an isolated point. 

\begin{corollary} \label{cor10}
Let $A\in B(\eH)$. If $\eM \in \Lat(A)$ is such that $P_{_\eM}$ commutes with every $P_{_\eN}$ $(\eN\in \Lat(A))$, then $\eM$ is an ultrainvariant
subspace of $A$.
\end{corollary}

\begin{proof}
Assume that $\eM \in \Lat(A)$ is such that $P_{_\eM}P_{_\eN}=P_{_\eN}P_{_\eM}$ for every $\eN\in \Lat(A)$. Then, for every $\eN\in \Lat(A)$, 
$P_{_\eM}P_{_\eN}$ is an orthogonal projection such that $P_{_\eM}P_{_\eN}\leq P_{_\eM}$ and $P_{_\eM}P_{_\eN}\leq P_{_\eN}$. Suppose that $\eN\ne \eM$.
Then either $P_{_\eM}P_{_\eN}< P_{_\eM}$ or $P_{_\eM}P_{_\eN}< P_{_\eN}$. If the former holds, then there exists $x\in \Image(P_{_\eM})$, $\| x\|=1$,
such that $P_{_\eM}P_{_\eN}x=0$. Since $P_{_\eM}x=x$ it follows that $P_{_\eN}x=P_{_\eN}P_{_\eM}x=0$. Thus, $\| P_{_\eM}-P_{_\eN}\|\geq
\| P_{_\eM}x-P_{_\eN}x\|=\| x\|=1$ which gives $\| P_{_\eM}-P_{_\eN}\|=1$. We get the same conclusion if $P_{_\eM}P_{_\eN}< P_{_\eN}$ holds. We conclude that
$\eM$ is an isolated point of $\Lat(A)$. Hence, it is inaccessible and therefore it is an ultrainvariant subspace of $A$, by Theorem \ref{cor10}.
\end{proof}

Recall that $A\in B(\eH)$ is unicellular if $\Lat(A)$ is totally ordered, that is, for any two subspaces $\eM, \eN \in \Lat(A)$, either $\eM \subseteq \eN$
or $\eM\supseteq \eN$.

\begin{corollary} \label{cor11}
If $A\in B(\eH)$ is unicellular, then every $\eM \in \Lat(A)$ is an ultrainvariant subspace of $A$.
\end{corollary}

\begin{proof}
Assume that $A\in B(\eH)$ is unicellular and let $\eM \in \Lat(A)$ be arbitrary. If $\eN \in \Lat(A)$, then either $\eM \subseteq \eN$
or $\eM\supseteq \eN$ which implies that $P_{_\eM}$ and $P_{_\eN}$ commute. By Corollary \ref{cor10}, $\eM$ is an ultrainvariant subspace of $A$.
\end{proof}

An example of a unicellular operator is the Volterra operator $V$ on $L^2[0,1]$. By \cite[Example 3]{Don} (see also \cite[Theorem 4.14]{RR}),
$\Lat(V)$ consists of subspaces $\eM_\alpha=\{f\in L^2[0,1];\; f=0\;\text{a.e. on}\; [0,\alpha]\}$, where $\alpha \in [0,1]$. Hence, by Corollary \ref{cor11},
every invariant subspace of $V$ is actually an ultrainvariant subspace.

Let $(e_n)_{n=0}^{\infty}$ be an orthonormal basis of $\eH$ and let $w=(w_n)_{n=1}^{\infty}$ be a bounded sequence of complex numbers. Then $w$ determines
an operator $W_w\in \eH$ by $W_w e_0=0$ and $W_w e_n=w_n e_{n-1}$ for $n\in \bN$ (hence, $W_w$ is the weighted backward shift determined by $w$).
For each $n\geq 0$, denote $\eM_n=\bigvee\{ e_0, \ldots, e_n\}$.

\begin{proposition} \label{prop11}
Every $\eM_n\subseteq \eH$ $(n\geq 0)$ is an ultrainvariant subspace of $W_w$.
\end{proposition}

\begin{proof}
Since $\eM_0=\Ker(W_w)$ and kernels are ultrainvariant (see the paragraph before Lemma \ref{lem03}) we may assume that $n\geq 1$. Observe that $W_w x\in \eM_k$ for 
some $x\in \eH$ and $k\geq 0$ if and only if $x\in \eM_{k+1}$. Indeed, if $x=\sum\limits_{j=0}^{\infty} (x,e_j)e_j$ and $(x,e_j)\ne 0$ for some $j>k+1$, then
$W_w x=\sum_{j=1}^{\infty} (x,e_j) w_j e_{j-1}\not\in \eM_{k}$.

Let $T\in \LCom(W_w;\eM_n)$ be arbitrary. Then $W_w Te_0=TW_w e_0=0$ and therefore $Te_0\in \eM_0$. Suppose that for $j<n$ we have $Te_j\in \eM_j$. Then
$W_w T e_{j+1}=TW_w e_{j+1}=w_jTe_j\in \eM_j$ which gives $Te_{j+1}\in \eM_{j+1}$. We conclude that $\eM_n$ is invariant for every $T\in \LCom(W_w;\eM_n)$,
that is, $\eM_n$ is ultrainvariant.
\end{proof}

If $w$ is a sequence of non-zero complex numbers such that $|w_1|>|w_2|>\cdots$ and $\sum\limits_{n=1}^{\infty}|w_n|^2<\infty$, then $W_w$
is called the Donoghue operator with weight sequence $w$. By \cite[Theorem 4.12]{RR} (a special case was considered by Donoghue 
\cite[Example 1]{Don}), a closed subspace $\eM\subseteq \eH$ is a non-trivial invariant subspace of the Donoghue operator $W_w$ if and only if
$\eM=\eM_n$ for some $n\geq 0$. Hence the Donoghue operator is unicellular and therefore every $\eM_n$ is an ultrainvariant subspace of $W_w$.\bigskip

{\bf Acknowledgments.} 
The author was supported by the Slovenian Research Agency through the research program P2-0268.

%%%%%%%%%%%%%%%%%%%%%%%%%%%%%%%%%%%%%%%%%%%%%%%%%%%%%%%%%%%%%%%%%%%%%%%%%%%%%%%%%%%%%%%%
%%%%%%%%%%%%%%%%%%%%%%%%%%%%%%%%%%%%%%%%%%%%%%%%%%%%%%%%%%%%%%%%%%%%%%%%%%%%%%%%%%%%%%%%

\end{document}